\font\tenmib=cmmib10
\font\sevenmib=cmmib10 scaled 800
\font\titolo=cmbx12
\font\titolone=cmbx10 scaled\magstep 2

\font\cs=cmcsc10
\font\sc=cmcsc10
\font\css=cmcsc8

\font\ninerm=cmr9
\font\ottorm=cmr8
\textfont5=\tenmib
\scriptfont5=\sevenmib
\scriptscriptfont5=\fivei

\font\euftw=eufm9 scaled\magstep1
\font\euftww=eufm7 scaled\magstep1
\font\euftwww=eufm5 scaled\magstep1
\font\msytw=msbm9 scaled\magstep1

\font\msytwww=msbm6 scaled\magstep1

\font\indbf=cmbx10 scaled\magstep2

\font\ottorm=cmr8\font\ottoi=cmmi8\font\ottosy=cmsy8%
\font\ottobf=cmbx8\font\ottott=cmtt8%
\font\ottocss=cmcsc8%
\font\ottosl=cmsl8\font\ottoit=cmti8%
\font\sixrm=cmr6\font\sixbf=cmbx6\font\sixi=cmmi6\font\sixsy=cmsy6%
\font\fiverm=cmr5\font\fivesy=cmsy5
\font\fivei=cmmi5
\font\fivebf=cmbx5%

\def\ottopunti{\def\rm{\fam0\ottorm}%
\textfont0=\ottorm\scriptfont0=\sixrm\scriptscriptfont0=\fiverm%
\textfont1=\ottoi\scriptfont1=\sixi\scriptscriptfont1=\fivei%
\textfont2=\ottosy\scriptfont2=\sixsy\scriptscriptfont2=\fivesy%
\textfont3=\tenex\scriptfont3=\tenex\scriptscriptfont3=\tenex%
\textfont4=\ottocss\scriptfont4=\sc\scriptscriptfont4=\sc%
\textfont5=\tenmib\scriptfont5=\sevenmib\scriptscriptfont5=\fivei
\textfont\itfam=\ottoit\def\it{\fam\itfam\ottoit}%
\textfont\slfam=\ottosl\def\sl{\fam\slfam\ottosl}%
\textfont\ttfam=\ottott\def\tt{\fam\ttfam\ottott}%
\textfont\bffam=\ottobf\scriptfont\bffam=\sixbf%
\scriptscriptfont\bffam=\fivebf\def\bf{\fam\bffam\ottobf}%
\setbox\strutbox=\hbox{\vrule height7pt depth2pt width0pt}%
\normalbaselineskip=9pt\let\sc=\sixrm\normalbaselines\rm}
\let\nota=\ottopunti%

%
%
%
%
%
%
%

\global\newcount\numsec\global\newcount\numapp
\global\newcount\numfor\global\newcount\numfig
\global\newcount\numsub
\numsec=0\numapp=0\numfig=0
\def\veroparagrafo{\number\numsec}\def\veraformula{\number\numfor}
\def\veraappendice{\number\numapp}\def\verasub{\number\numsub}
\def\verafigura{\number\numfig}

\def\section(#1,#2){\advance\numsec by 1\numfor=1\numsub=1\numfig=1%
\SIA p,#1,{\veroparagrafo} %
\write15{\string\Fp (#1){\secc(#1)}}%
\write16{ sec. #1 ==> \secc(#1)  }%
\hbox to \hsize{\titolo\hfill \number\numsec. #2\hfill%
\expandafter{\alato(sec. #1)}}\*}

\def\appendix(#1,#2){\advance\numapp by 1\numfor=1\numsub=1\numfig=1%
\SIA p,#1,{A\veraappendice} %
\write15{\string\Fp (#1){\secc(#1)}}%
\write16{ app. #1 ==> \secc(#1)  }%
\hbox to \hsize{\titolo\hfill Appendix A\number\numapp. #2\hfill%
\expandafter{\alato(app. #1)}}\*}

\def\senondefinito#1{\expandafter\ifx\csname#1\endcsname\relax}

\def\SIA #1,#2,#3 {\senondefinito{#1#2}%
\expandafter\xdef\csname #1#2\endcsname{#3}\else
\write16{???? ma #1#2 e' gia' stato definito !!!!} \fi}

\def \Fe(#1)#2{\SIA fe,#1,#2 }
\def \Fp(#1)#2{\SIA fp,#1,#2 }
\def \Fg(#1)#2{\SIA fg,#1,#2 }

\def\etichetta(#1){(\veroparagrafo.\veraformula)%
\SIA e,#1,(\veroparagrafo.\veraformula) %
\global\advance\numfor by 1%
\write15{\string\Fe (#1){\equ(#1)}}%
\write16{ EQ #1 ==> \equ(#1)  }}

\def\etichettaa(#1){(A\veraappendice.\veraformula)%
\SIA e,#1,(A\veraappendice.\veraformula) %
\global\advance\numfor by 1%
\write15{\string\Fe (#1){\equ(#1)}}%
\write16{ EQ #1 ==> \equ(#1) }}

\def\getichetta(#1){\veroparagrafo.\verafigura%
\SIA g,#1,{\veroparagrafo.\verafigura} %
\global\advance\numfig by 1%
\write15{\string\Fg (#1){\graf(#1)}}%
\write16{ Fig. #1 ==> \graf(#1) }}

\def\etichettap(#1){\veroparagrafo.\verasub%
\SIA p,#1,{\veroparagrafo.\verasub} %
\global\advance\numsub by 1%
\write15{\string\Fp (#1){\secc(#1)}}%
\write16{ par #1 ==> \secc(#1)  }}

\def\etichettapa(#1){A\veraappendice.\verasub%
\SIA p,#1,{A\veraappendice.\verasub} %
\global\advance\numsub by 1%
\write15{\string\Fp (#1){\secc(#1)}}%
\write16{ par #1 ==> \secc(#1)  }}

\def\Eq(#1){\eqno{\etichetta(#1)\alato(#1)}}
\def\eq(#1){\etichetta(#1)\alato(#1)}
\def\Eqa(#1){\eqno{\etichettaa(#1)\alato(#1)}}
\def\eqa(#1){\etichettaa(#1)\alato(#1)}
\def\eqg(#1){\getichetta(#1)\alato(fig. #1)}
\def\sub(#1){\0\palato(p. #1){\bf \etichettap(#1).}}
\def\asub(#1){\0\palato(p. #1){\bf \etichettapa(#1).}}

\def\equv(#1){\senondefinito{fe#1}$\clubsuit$#1%
\write16{eq. #1 non e' (ancora) definita}%
\else\csname fe#1\endcsname\fi}
\def\grafv(#1){\senondefinito{fg#1}$\clubsuit$#1%
\write16{fig. #1 non e' (ancora) definito}%
\else\csname fg#1\endcsname\fi}
\def\secv(#1){\senondefinito{fp#1}$\clubsuit$#1%
\write16{par. #1 non e' (ancora) definito}%
\else\csname fp#1\endcsname\fi}

\def\equ(#1){\senondefinito{e#1}\equv(#1)\else\csname e#1\endcsname\fi}
\def\graf(#1){\senondefinito{g#1}\grafv(#1)\else\csname g#1\endcsname\fi}
\def\figura(#1){{\css Figura} \getichetta(#1)}
\def\secc(#1){\senondefinito{p#1}\secv(#1)\else\csname p#1\endcsname\fi}
\def\sec(#1){{\S\secc(#1)}}
\def\refe(#1){{[\secc(#1)]}}

\def\BOZZA{\bz=1
\def\alato(##1){\rlap{\kern-\hsize\kern-1.2truecm{$\scriptstyle##1$}}}
\def\palato(##1){\rlap{\kern-1.2truecm{$\scriptstyle##1$}}}
				  }

\def\alato(#1){}
\def\galato(#1){}
\def\palato(#1){}


{\count255=\time\divide\count255 by 60 \xdef\hourmin{\number\count255}
        \multiply\count255 by-60\advance\count255 by\time
   \xdef\hourmin{\hourmin:\ifnum\count255<10 0\fi\the\count255}}

\def\oramin{\hourmin }

\def\data{\number\day/\ifcase\month\or gennaio \or febbraio \or marzo \or
aprile \or maggio \or giugno \or luglio \or agosto \or settembre
\or ottobre \or novembre \or dicembre \fi/\number\year;\ \oramin}

\newdimen\xshift \newdimen\xwidth \newdimen\yshift \newdimen\ywidth

\def\ins#1#2#3{\vbox to0pt{\kern-#2\hbox{\kern#1 #3}\vss}\nointerlineskip}

\def\eqfig#1#2#3#4#5{
\par\xwidth=#1 \xshift=\hsize \advance\xshift
by-\xwidth \divide\xshift by 2
\yshift=#2 \divide\yshift by 2
\line{\hglue\xshift \vbox to #2{\vfil
#3 \includegraphics{#4.ps}
}\hfill\raise\yshift\hbox{#5}}}

\def\8{\write12}  


\let\a=\alpha \let\b=\beta  \let\g=\gamma  \let\d=\delta \let\e=\varepsilon
  \let\h=\eta   \let\th=\theta \let\k=\kappa \let\l=\lambda
    \let\n=\nu             
 \let\t=\tau    
   \let\o=\omega
\let\G=\Gamma \let\D=\Delta  \let\Th=\Theta\let\L=\Lambda \let\X=\Xi
         
\let\O=\Omega 

\def\\{\hfill\break} \let\==\equiv

\let\io=\infty 

\let\0=\noindent

\let\dpr=\partial

\def\tende#1{\,\vtop{\ialign{##\crcr\rightarrowfill\crcr
 \noalign{\kern-1pt\nointerlineskip} \hskip3.pt${\scriptstyle
 #1}$\hskip3.pt\crcr}}\,}
\def\circage{\lower2pt\hbox{$\,\buildrel > \over {\scriptstyle \sim}\,$}}
\def\otto{\,{\kern-1.truept\leftarrow\kern-5.truept\to\kern-1.truept}\,}

\def\EE{{\cal E}}\def\MM{{\cal M}} \def\VV{{\cal V}}
\def\CC{{\cal C}}\def\FF{{\cal F}} \def\HHH{{\cal H}}
\def\NN{{\cal N}} \def\BBB{{\cal B}}\def\II{{\cal I}}
\def\RR{{\cal R}}  
\def\AAA{{\cal A}} \def\SS{{\cal S}}

\def\T#1{{#1_{\kern-3pt\lower7pt\hbox{$\widetilde{}$}}\kern3pt}}
\def\VVV#1{{\underline #1}_{\kern-3pt
\lower7pt\hbox{$\widetilde{}$}}\kern3pt\,}
\def\W#1{#1_{\kern-3pt\lower7.5pt\hbox{$\widetilde{}$}}\kern2pt\,}

\def\lis{\overline}

\def\indica{\leaders \hbox to 0.5cm{\hss.\hss}\hfill}
\def\guida{\leaders\hbox to 1em{\hss.\hss}\hfill}

\def\hh{{\bf h}}  \def\AA{{\bf A}} 
\def\BB{{\bf B}}   
   
\def\aaa{{\bf a}}\def\bbb{{\bf b}}\def\III{{\bf I}}

\def\ul{\underline}


\mathchardef\aa   = "050B
\mathchardef\bb   = "050C
\mathchardef\ggg  = "050D
\mathchardef\xxx  = "0518
\mathchardef\zzzzz= "0510
\mathchardef\oo   = "0521
\mathchardef\lll  = "0515
\mathchardef\mm   = "0516
\mathchardef\Dp   = "0540
\mathchardef\H    = "0548
\mathchardef\FFF  = "0546
\mathchardef\ppp  = "0570
\mathchardef\nn   = "0517
\mathchardef\ff   = "0527
\mathchardef\pps  = "0520
\mathchardef\FFF  = "0508
\mathchardef\nnnnn= "056E

\def\to{\rightarrow}

\def\qed{\hfill\raise1pt\hbox{\vrule height5pt width5pt depth0pt}}

\def\Val{{\rm Val}}
\def\indic{\hbox{\raise-2pt \hbox{\indbf 1}}}

\def\RRR{\hbox{\msytw R}}

\def\NNN{\hbox{\msytw N}} 
 \def\ZZZ{\hbox{\msytw Z}}
 \def\zzz{\hbox{\msytwww Z}}
\def\TTT{\hbox{\msytw T}}

    \def\vvvv{\hbox{\euftww v}}
\def\vvvvv{\hbox{\euftwww v}}
  \def\wwwww{\hbox{\euftwww w}}

\def\MMMM{\hbox{\msytw M}}
\def\AAAA{\hbox{\msytw A}}
\def\BBBB{\hbox{\euftw B}}

\def\ul#1{{\underline#1}}

\def\V0{{\bf 0}}


\newcount\mgnf  
\mgnf=0 

\ifnum\mgnf=0
\def\openone{\leavevmode\hbox{\ninerm 1\kern-3.3pt\tenrm1}}%
\def\*{\vglue0.3truecm}\fi
\ifnum\mgnf=1
\def\openone{\leavevmode\hbox{\ninerm 1\kern-3.63pt\tenrm1}}%
\def\*{\vglue0.5truecm}\fi


\newcount\tipobib\newcount\bz\bz=0\newcount\aux\aux=1
\newdimen\bibskip\newdimen\maxit\maxit=0pt


\tipobib=0
\def\9#1{\ifnum\aux=1#1\else\relax\fi}

\newwrite\bib
\immediate\openout\bib=\jobname.bib
\global\newcount\bibex
\bibex=0
\def\verabib{\number\bibex}

\ifnum\tipobib=0
\def\cita#1{\expandafter\ifx\csname c#1\endcsname\relax
\hbox{$\clubsuit$}#1\write16{Manca #1 !}%
\else\csname c#1\endcsname\fi}
\def\rife#1#2#3{\immediate\write\bib{\string\raf{#2}{#3}{#1}}
\immediate\write15{\string\C(#1){[#2]}}
\setbox199=\hbox{#2}\ifnum\maxit < \wd199 \maxit=\wd199\fi}
\fi
\ifnum\tipobib=1
\def\cita#1{%
\expandafter\ifx\csname d#1\endcsname\relax%
\expandafter\ifx\csname c#1\endcsname\relax%
\hbox{$\clubsuit$}#1\write16{Manca #1 !}%
\else\probib(ref. numero )(#1)%
\csname c#1\endcsname%
\fi\else\csname d#1\endcsname\fi}%
\def\rife#1#2#3{\immediate\write15{\string\Cp(#1){%
\string\immediate\string\write\string\bib{\string\string\string\raf%
{\string\verabib}{#3}{#1}}%
\string\Cn(#1){[\string\verabib]}%
\string\CCc(#1)%
}}}%
\fi
\ifnum\tipobib=2%
\def\cita#1{\expandafter\ifx\csname c#1\endcsname\relax
\hbox{$\clubsuit$}#1\write16{Manca #1 !}%
\else\csname c#1\endcsname\fi}
\def\rife#1#2#3{\immediate\write\bib{\string\raf{#1}{#3}{#2}}
\immediate\write15{\string\C(#1){[#1]}}
\setbox199=\hbox{#2}\ifnum\maxit < \wd199 \maxit=\wd199\fi}
\fi

\def\Cn(#1)#2{\expandafter\xdef\csname d#1\endcsname{#2}}
\def\CCc(#1){\csname d#1\endcsname}
\def\probib(#1)(#2){\global\advance\bibex+1%
\9{\immediate\write16{#1\verabib => #2}}%
}

\def\C(#1)#2{\SIA c,#1,{#2}}
\def\Cp(#1)#2{\SIAnx c,#1,{#2}}

\def\SIAnx #1,#2,#3 {\senondefinito{#1#2}%
\expandafter\def\csname#1#2\endcsname{#3}\else%
\write16{???? ma #1,#2 e' gia' stato definito !!!!}\fi}

\bibskip=10truept
\def\hboxto{\hbox to}

\catcode`\{=12\catcode`\}=12
\catcode`\<=1\catcode`\>=2
\immediate\write\bib<
        \string\halign{\string\hboxto \string\maxit%
        {\string #\string\hfill}&%
        \string\vtop{\string\parindent=0pt\string\advance\string\hsize%
        by -.5truecm%
        \string#\string\vskip \bibskip
        }\string\cr%
>
\catcode`\{=1\catcode`\}=2
\catcode`\<=12\catcode`\>=12

\def\raf#1#2#3{\ifnum \bz=0 [#1]&#2 \cr\else
\llap{${}_{\rm #3}$}[#1]&#2\cr\fi}

\newread\bibin

\catcode`\{=12\catcode`\}=12
\catcode`\<=1\catcode`\>=2
\def\chiudibib<
\catcode`\{=12\catcode`\}=12
\catcode`\<=1\catcode`\>=2
\immediate\write\bib<}>
\catcode`\{=1\catcode`\}=2
\catcode`\<=12\catcode`\>=12
>
\catcode`\{=1\catcode`\}=2
\catcode`\<=12\catcode`\>=12

\def\makebiblio{
\ifnum\tipobib=0
\advance \maxit by 10pt
\else
\maxit=1.truecm
\fi
\chiudibib
\immediate \closeout\bib
\openin\bibin=\jobname.bib
\ifeof\bibin\relax\else
\raggedbottom
\input \jobname.bib
\fi
}

\openin13=#1.aux \ifeof13 \relax \else
\input #1.aux \closein13\fi
\openin14=\jobname.aux \ifeof14 \relax \else
\input \jobname.aux \closein14 \fi
\immediate\openout15=\jobname.aux

\def\biblio{\*\*\centerline{\titolo References}\*\nobreak\makebiblio}


\ifnum\mgnf=0
   \magnification=\magstep0
   \hsize=14.truecm\vsize=18.0truecm\voffset2.truecm\hoffset.5truecm
   \parindent=0.3cm\baselineskip=0.45cm\fi
\ifnum\mgnf=1
   \magnification=\magstep1\hoffset=0.truecm
   \hsize=14.truecm\vsize=18.0truecm
   \baselineskip=18truept plus0.1pt minus0.1pt \parindent=0.9truecm
   \lineskip=0.5truecm\lineskiplimit=0.1pt      \parskip=0.1pt plus1pt\fi


\mgnf=0   
\openin14=\jobname.aux \ifeof14 \relax \else
\input \jobname.aux \closein14 \fi
\openout15=\jobname.aux

\footline={\rlap{\hbox{\copy200}}\tenrm\hss \number\pageno\hss}
\def\fiat{}

%
\fiat




\centerline{\titolone Degenerate lower-dimensional tori}
\vskip.1truecm
\centerline{\titolone under the Bryuno condition}

\vskip.5truecm
\centerline{{\titolo Guido Gentile}}
\vskip.2truecm
\centerline{Dipartimento di Matematica,
Universit\`a di Roma Tre, Roma, I-00146, Italy}
\vskip1.0truecm

\line{\vtop{
\line{\hskip1truecm\vbox{\advance \hsize by -2.1 truecm
\0{\cs Abstract.}
{\it We study the problem of conservation of maximal and lower-dimensional
invariant tori for analytic convex quasi-integrable Hamiltonian systems.
In the absence of perturbation the lower-dimensional tori are degenerate,
in the sense that the normal frequencies vanish, so that the tori are
neither elliptic nor hyperbolic. We show that if the perturbation
parameter is small enough, for a large measure subset of any resonant
submanifold of the action variable space, under some generic
non-degeneracy conditions on the perturbation function,
there are lower-dimensional tori which are conserved.
They are characterised by rotation vectors satisfying some generalised
Bryuno conditions involving also the normal frequencies.
We also show that, again under some generic assumptions on the
perturbation, any torus with fixed rotation vector satisfying
the Bryuno condition is conserved for most values of the perturbation
parameter in an interval small enough around the origin. 
According to the sign of the normal frequencies and of the
perturbation parameter the torus
becomes either hyperbolic or elliptic or of mixed type.}} \hfill}
}}

\*\*\*\*
\section(1,Introduction)

\0It is well known that in quasi-integrable analytic Hamiltonian systems
KAM invariant tori are conserved under conditions on the
rotation vectors milder than the usual Diophantine condition
originally introduced by Kolmogorov \cita{Ko}.
A more general condition was introduced by Bryuno in Refs.~\cita{Br1}
and \cita{Br2}, and it is nowadays knowns as the {\it Bryuno condition}.
Among the most exhaustive studies in this direction
we cite those by R\"ussmann (for a recent review see Ref.~\cita{R2}).
In some related problems, such as Siegel's problem
(in the analytic framework), one knows that the Bryuno condition is a
necessary and sufficient condition for the dynamics to be conjugated
to the linear one, as the work by Yoccoz has shown \cita{Y}.
In the case of area-preserving maps the same result likely holds,
and for the standard map this has been explicitly verified.
For Siegel's problem, Yoccoz also proved that a deep relationship exists
between the radius of convergence of the linearising function and the
so-called Bryuno function. An analogous relationship between the
radius of convergence of the conjugating function and the Bryuno function
has been found for the standard map by combining the results
of Davie \cita{D} with those of Berretti and Gentile \cita{BeG}.
We mention also the work by Ecalle and Vallet \cita{EV}, where it is shown
that, under the Bryuno condition, all analytic resonant vector fields
and diffeomorphisms admit an analytic correction which make them
linearisable (as conjectured by Gallavotti \cita{Ga1}, and proved under
the usual Diophantine condition by Eliasson \cita{E2}, hence by
Gentile and Mastropietro \cita{GM} with techniques more similar to those
we use in the present paper). Note that for such a problem
the rotation vector is fixed  and no value of the perturbation parameter
has to be excluded, so that the problem rather simplifies, as all
the difficulties related to estimating the measure of the
allowed values for the parameters disappear.
For instance there is not the difficulty of including
the correction in the original vector field with a different
unperturbed rotation vector, as in the case of the KAM theorem
for isochronous systems (cf. Ref.~\cita{BaG} for a discussion
within the formalism used here). Extensions of the Bryuno condition
to other contexts, such as the reducibility of skew-products,
has been partially provided recently by Lopes Dias \cita{LD}.

In this paper we consider a problem of lower-dimensional tori
similar to that considered in Refs.~\cita{P2} and \cita{R2},
with the main difference being that the normal frequencies vanish
in the absence of perturbation. Such a problem has been
explicitly considered in a series of papers, such as
Refs. \cita{T}, \cita{JLZ}, \cita{Ch1}, \cita{GG1} and \cita{GG2}.
We refer to the latter for an introduction,
and for a review of the existing results.
In particular we start by considering the same class of Hamiltonian
systems with $d\ge 2$ degrees of freedom considered in Refs.~\cita{GG1}
and \cita{GG2}, originally introduced in Ref.~\cita{Th},
$$ \HHH = {1\over2}\AA\cdot\AA +
{1 \over 2} \BB\cdot\BB + \e f(\aa,\bb) ,
\Eq(1.1) $$
where $(\aa,\AA)\in \TTT^{r}\times \RRR^{r}$ and
$(\bb,\BB)\in \TTT^{s}\times \RRR^{s}$ are conjugate action-angle
variables, with $r+s=d$, and $\cdot$ denotes the inner product
both in $\RRR^{r}$ and in $\RRR^{s}$.
The perturbation $f(\aa,\bb)$ is assumed to be real analytic,
so that, if we write
$$ f(\aa,\bb) = \sum_{\nn\in\zzz^{r}} {\rm e}^{i\nn\cdot\aa}
f_{\nn}(\bb) ,
\Eq(1.2) $$
there exist positive constants $F_{0}$, $F_{1}$ and $\k_{0}$
such that $|\dpr_{\bb}^{q}f_{\nn}(\bb)| \le
q!F_{0}F_{1}^{q}{\rm e}^{-\k_{0} |\nn|}$ for all $\nn\in\ZZZ^{r}$,
all $\bb\in\TTT^{s}$ and all $q\in\ZZZ_{+}$.
For $\bb_{0}$ a stationary point of $f_{\V0}(\bb)$ we call
$a_{1},\ldots,a_{s}$ the eigenvalues of the
matrix $\dpr_{\bb}^{2}f_{\V0}(\bb_{0})$.
The case of maximal tori is recovered by setting $r=d$.
The general case of Hamiltonians describing perturbations
of any convex systems will be briefly discussed
in Appendix \secc(A2), even if the full discussion is
deferred to Ref.~\cita{GGG2}. Here we prefer to concentrate ourselves
to the simpler model \equ(1.1), in order to distinguish between
the more relevant features of the renormalisation group techniques
and the more technical intricacies pertaining rather to problems
of spectral analysis and matrix algebra.

For $\e=0$ the system described by the Hamiltonian \equ(1.1)
is integrable. Any solution of the form
$(\aa,\bb,\AA,\BB)=(\aa_{0}+\oo t,\bb_{0},\AA_{0},\BB_{0})$,
with $\oo=\AA_{0}$ having rationally independent components,
fills densely a lower-dimensional torus with rotation vector $\oo$.
We call $\AA$ and $\BB$ the {\it non-resonant} and {\it resonant},
respectively, action variables. With a shift of the
resonant action variables we can always assume $\BB_{0}=\V0$.
The so-called {\it normal frequencies}, that is the frequencies
describing the dynamics of the $(\bb,\BB)$-variables, vanish for $\e=0$,
so that the considered unperturbed torus is neither elliptic
nor hyperbolic (nor of mixed type). We refer to such a situation
by saying that one has a {\it degenerate torus}.

The frequency map $\AA \mapsto \oo(\AA)$ is a local diffeomorphism
(in our case it is trivially the identity), so that the condition
$\BB=\V0$ defines an $r$-dimensional manifold $\MMMM_{r}$
({\it resonant submanifold}), which is determined by the space
of the non-resonant action variable $\AA$; we call the latter
the {\it non-resonant action variable space}.

We are interested in two different problems.\\
(A) One can fix the perturbation parameter (small enough) and
study for which rotation vectors some invariant tori are conserved,
in the spirit of the KAM theorem for maximal tori,
and as done in most of the papers on such a subject, as
Refs. \cita{M1}, \cita{M2}, \cita{Ku1}, \cita{Ku2}, \cita{E1}, \cita{P2},
\cita{Ch2}, \cita{R2}, \cita{BKS}, \cita{JV2}, and many others.\\
(B) Either one can look at a lower-dimensional invariant torus
with fixed rotation vector, and study the dependence of such a torus
on the perturbation parameter. For instance this has been done,
with the techniques used here, in Refs.~\cita{GG1} and \cita{GG2}.

The same twofold program has been followed, under the usual Diophantine
condition, in Refs.~\cita{Ge2} and \cita{Ge3} in the study of
the quasi-periodic solutions and of the spectrum for a class of
two-level systems in a strong quasi-periodic external field.

About problem (A) we find the analogue of R\"ussmann's \cita{R2}
and P\"oschel's \cita{P2} result for systems with
distinct\footnote{${}^1$}{\nota We prefer using the word `distinct',
instead of `non-degenerate', just to avoid confusion, as we are calling
the normal frequencies `degenerate' when they vanish for $\e=0$.}
normal frequencies of order 1.
In addition, our result applies also
to the case of non-distinct normal frequencies,
provided that they are all different from zero and of order $\e$,
that is provided degeneracy is removed to first order. In particular
this means that our result does not follow from the works
available in literature: in principle one could think to perform a
canonical transformation which introduces normal frequencies
of order $\sqrt{|\e|}$, while keeping the perturbation
to order $\e$ (as explicitly done in Refs.~\cita{T} and \cita{ChW}),
but in this way the normal frequencies are still required
to be distinct in order to apply P\"oschel's result,
whereas we do not need such a condition.
On the other hand the case of possibly non-distinct normal frequencies
(of order 1) has been dealt with in Refs.~\cita{Bo1},
\cita{Bo2}, \cita{You}, \cita{Xu} and \cita{XY}
only under the usual Diophantine condition.

\*

\0{\bf Theorem 1.} {\it Consider the Hamiltonian \equ(1.1).
Suppose $\bb_{0}$ to be such that $\dpr_{\bb}f_{\V0}(\bb_{0})=0$,
and assume that the eigenvalues of the matrix
$\dpr_{\bb}^{2}f_{\V0}(\bb_{0})$ are all different from zero
(that is $a_{i}\neq0$ for all $i=1,\ldots,s$).
Let $\AAA \subset \RRR^{r}$ be any open set of the
non-resonant action variable space. Then for any $\d>0$
there are $\e_{0}$ small enough and a subset $\AAA_{*}\subset\AAA$
such that if $|\e|<\e_{0}$ and $\oo\in\AAA_{*}$
the system described by the Hamiltonian \equ(1.1)
admits a lower-dimensional torus of the form
$$ \cases{
\aa = \pps + \aaa (\pps,\bb_{0},\oo,\e) , \cr
\bb = \bb_{0} + \bbb (\pps,\bb_{0},\oo,\e) , \cr
\AA = \oo + \left( \oo \cdot \dpr_{\pps} \right)
\aaa (\pps,\bb_{0},\oo,\e) , \cr
\BB = \left( \oo \cdot \dpr_{\pps} \right)
\bbb (\pps,\bb_{0},\oo,\e) , \cr}
\Eq(1.3) $$
with the functions $\aaa$ and $\bbb$ vanishing at $\e=0$,
analytic and periodic in $\pps$, and the Lebesgue measure of the set
$\AAA\setminus\AAA^{*}$ is less than $\d$.
The parameterisation in \equ(1.3) is such that
$\pps=\pps_{0}+\oo t$ describes a linear flow on $\TTT^{r}$.
In the case of maximal tori ($r=d$) the subset of phase space
which is filled by invariant tori has complement whose Lebesgue measure
is less than $C\d$, for some positive constant $C$.}

\*

We can interpret Theorem 1 by saying that, in the presence of
perturbations, the resonant tori are destroyed in general,
but some of them survive. They are determined by the stationary
points of the {\it potential function} $f_{\V0}(\bb)$. Let
$\bb_{0}$ one of these stationary points: under the (generic)
non-degeneracy conditions assumed on the eigenvalues
of the matrix $\dpr_{\bb}^{2}f_{\V0}(\bb_{0})$, we can say that,
in correspondence of such a point $\bb_{0}$,
there is a conserved invariant lower-dimensional torus.
The latter is either {\it hyperbolic} or {\it elliptic} or
{\it of mixed type} according to the signs of the eigenvalues and of the
perturbation parameter $\e$: it is elliptic if the eigenvalues
have the same sign as $\e$, hyperbolic if they have opposite sign
with respect to $\e$, and of mixed type otherwise. The rotation vector
of such a torus will be found to satisfy some Diophantine conditions
involving also the normal frequencies. In particular we shall find that
any maximal torus with rotation vector $\oo$ which is a Bryuno vector
in $\RRR^{d}$ is conserved provided that $\e$ is small enough
(depending on $\oo$). By {\it Bryuno vector} in $\RRR^{r}$
we mean a vector $\oo\in\RRR^{r}$ such that
$$ \sum_{n=0}^{\io} {1 \over 2^{n}}
\log {1 \over
{\displaystyle \inf_{0<|\nn|\le 2^{n}} |\oo\cdot\nn|} } < \io .
\Eq(1.4) $$
For further properties of Bryuno vectors we refer to Section \secc(2).
Note that in Refs.~\cita{P2} and \cita{R2} slightly different (but
equivalent) conditions are found for the rotation vectors of the
surviving elliptic tori, which can be expressed in terms of a suitable
approximation function, first introduced by R\"ussmann~\cita{R1};
cf. the quoted references for further details.

Also concerning problem (B) our result does not exist in literature,
and it represents the natural extension of Ref.~\cita{GG1}
and \cita{GG2} to the case of more general rotation vectors.
In this case we still need the condition for the normal frequencies
to be distinct, as in Ref.~\cita{GG2}.
Such a condition could be weakened by assuming only that degeneracy
is removed to some finite perturbation order; cf. Ref.~\cita{GGG1}.
Though, we shall impose on the rotation vector
only the Bryuno condition \equ(1.4),
a condition much weaker than Kolmogorov's Diophantine condition
for maximal tori and Mel$'$nikov condition for elliptic
lower-dimensional tori, as usually assumed
(cf. Refs.~\cita{T}, \cita{JLZ} and \cita{JV1}).

\*

\0{\bf Theorem 2.} {\it Let $\oo$ be a vector in $\RRR^{r}$
satisfying the Diophantine condition \equ(1.4).
Suppose $\bb_{0}$ to be such that $\dpr_{\bb}f_{\V0}(\bb_{0})=0$,
and assume that the eigenvalues of the matrix
$\dpr_{\bb}^{2}f_{\V0}(\bb_{0})$ are all different from zero
and pairwise distinct (that is $a_{i} \neq 0$ for all $i=1,\ldots,s$
and $a_{i} \neq a_{j}$ for all $1 \le i < j \le s$).
Then there exists $\e_{0}$ and a set $\EE\subset(-\e_{0},\e_{0})$,
with a density point at the origin,
such that for all $\e\in\EE$ there is a lower-dimensional
torus for the system described by the Hamiltonian \equ(1.1)
with rotation vector $\oo$, which can be parameterised as \equ(1.3),
with $\pps\in\TTT^{d}$ and the functions $\aaa$ and $\bbb$
vanishing at $\e=0$, analytic and periodic in $\pps$.}

\*

A density point for $\EE$ at the origin means that
the relative Lebesgue measure of the set $\EE\cap(-\e,\e)$,
that is ${\rm meas}(\EE\cap(-\e,\e))/2\e$, tends to 1 as $\e\to0$.
We shall say also, in such a case, that the set $\EE$ has
large relative (Lebesgue) measure in $(-\e_{0},\e_{0})$.

If $\e>0$ we can require $a_{i} \neq a_{j}$ for $i\neq j$
only for positive eigenvalues. If all eigenvalues are positive
then the corresponding torus is elliptic; if all eigenvalues
are negative then the corresponding torus is hyperbolic.
In the first case the allowed values of $\e$ form a Cantor
set in $[0,\e_{0})$ with large relative measure,
in the second one all values in $[0,\e_{0})$ are allowed.
The obvious analogue holds for $\e<0$.

For both theorems we shall give the proof in the case in which
all eigenvalues of $\dpr_{\bb}^{2}f_{\V0}(\bb_{0})$
are strictly positive and $\e>0$, which is the difficult case.
All the other cases can be obtained with trivial adaptations
of the proof. Note also that the case of maximal tori
can be obtained as a byproduct by setting $r=d$ in the following.

Finally we mention that if do not require that degeneracy
of the normal frequencies be removed to first order
then the problem can become much harder. For instance
if no condition at all is imposed on the perturbation
only partial results exist, and only for the case $s=1$
and $\oo$ a Diophantine rotation vector;
cf. Refs.~\cita{Ch1} and \cita{Ch2} (see also Ref.~\cita{GGG1}).

The paper is organised as follows. In Section \secc(2) we
introduce the Bryuno vectors, and we briefly review some properties
of theirs, which will be used in the forthcoming analysis.
Then Sections \secc(2) and \secc(3) are devoted
to the proof of Theorem 1 and of Theorem 2, respectively.
The proofs heavenly rely, both for notations and results,
on Ref.~\cita{GG2}, and we confine ourselves to give full details
only for the parts which are really different. In particular
the more technical aspects are deferred to Appendix \secc(A1).
However, by assuming the results of Ref.~\cita{GG2},
the discussion below becomes completely self-contained.
Finally in Appendix \secc(A2) we briefly discuss how the analysis
can be adapted to deal with more general convex Hamiltonian systems.

\*\*
\section(2,The Bryuno condition)

\0Given $\oo\in\RRR^{2}$ set $\o\=\min\{|\o_{1}|,|\o_{2}|\}/
\max\{|\o_{1}|,|\o_{2}|\}$. Let $\{q_{n}\}_{n=0}^{\io}$ be the
denominators of the convergents of $\o$.

The {\it Bryuno function} $\BBB(\o)$ is defined
as the solution of the functional equation \cita{Y}
$$ \cases{
\BBB(\o+1) = \BBB(\o) , & \cr
\BBB(\o)=-\log\o + \o \, \BBB(1/\o) , &
if $\o\in(0,1)$. \cr}
\Eq(2.1) $$
Define
$$ D(\oo) \= \sum_{n=0}^{\io} {\log q_{n+1}  \over q_{n}} .
\Eq(2.2) $$
Then it is easy to show that $\BBB(\o)<\io$ if and only if
$D(\oo)<\io$ \cita{Y}.

Given $\oo\in\RRR^{r}$ and $n\in\ZZZ_{+}$ set
$$ \a_{n}(\oo) = \inf_{0<|\nn|\le 2^{n}} \left| \oo\cdot\nn \right| ,
\Eq(2.3) $$
and define the {\it generalised Bryuno function} as
$$ B(\oo) = \sum_{n=0}^{\io} {1 \over 2^{n}}
\log {1 \over \a_{n}(\oo)} .
\Eq(2.4) $$

\*

\0{\bf Definition 1.}
{\it We shall call $\BBBB_{r} = \left\{ \oo \in \RRR^{r} : B(\oo) <
\io \right\}$ the set of {\it Bryuno vectors} in $\RRR^{r}$.
For any open set $\O\subset \RRR^{r}$ we call $\BBBB_{r}(\O)$ the
set of Bryuno vectors in $\O$.}

\*

The reason for this terminology relies on the following result.

\*

\0{\bf Lemma 1.} {\it For $r=2$ one has $\oo\in\BBBB_{2}$ if
and only if $D(\oo)<\io$.}

\*

\0{\it Proof.} Given $\oo\in\RRR^{2}$ assume for notation
simplicity $0<\o_{2}<\o_{1}$. Call $\a=\o_{2}/\o_{1}$,
and set $\oo=\o_{1}\oo_{0}$, with $\oo_{0}=(1,\a)$, so that $0<\a<1$
and $\log\a_{n}(\oo) = \log \o_{1} + \log \a_{n}(\oo_{0})$.
Consider the sequence of convergents $\{p_{n}/q_{n}\}_{n=1}^{\io}$
for $\a$ \cita{S}; one has $1/2q_{n+1}<|\a q_{n}-p_{n}|<1/q_{n+1}$,
and $|\oo_{0}\cdot\nn| > |\a q_{n}-p_{n}|$
for all $|\n_{2}|<q_{n+1}$, hence for all $|\nn|<2q_{n+1}$.

For each $n\ge 0$ define $r_{n}$ and $s_{n}$ such that
$2^{r_{n}-1} < 2 q_{n} \le 2^{r_{n}}$ and $r_{n}+s_{n}+1=r_{n+1}$.
Hence for all $r_{n} \le r' \le r_{n}+s_{n}$ one has
$\a_{r'}(\oo_{0}) = |\a q_{n}-p_{n}|$, which implies
$$  {1 \over 4} \left( {\log q_{n+1} \over q_{n}} \right) \le
\sum_{r'=r_{n}}^{r_{n}+s_{n}} {1 \over 2^{r'}}
\log {1 \over \a_{r'}(\oo_{0})} \le
{\log 2 \over q_{n}} + {\log q_{n+1} \over q_{n}} ,
\Eq(2.5) $$
so that, by using that $\sum_{n=0}^{\io} q_{n}^{-1} < \io$,
one obtains that there exist two positive constants $C_{1}$
and $C_{2}$ such that $D(\oo)/4 - C_{1} \le B(\oo) \le D(\oo) + C_{2}$,
and the assertion follows. \qed

\*

The sequence $\{\a_{n}(\oo)\}_{n=1}^{\io}$ is non-increasing,
so that it converges to 0 monotonically as $n\to\io$.
By taking possibly a subsequence, we can always
suppose $\a_{n+1}(\oo)<\a_{n}(\oo)$, strictly.

\*

\0{\bf Definition 2.}
{\it Set $\ZZZ^{r}_{*}=\ZZZ^{r}\setminus\{\V0\}$, and define
$$ n(\nn) = \left\{ n \in \ZZZ_{+} : 2^{n-1} <
|\nn| \le 2^{n} \right\}
= \inf \left\{ n \in \NNN : |\nn| \le 2^{n} \right\}
\Eq(2.6) $$
for any $\nn\in\ZZZ^{r}_{*}$.}

\*

For all $\nn\in\ZZZ^{r}_{*}$ one has, by definition,
$|\oo \cdot \nn| \ge \a_{n(\nn)}(\oo)$
and $2^{n(\nn)-1}<|\nn| \le 2^{n(\nn)}$.

\*

\0{\bf Definition 3.}
{\it Given a non-increasing sequence $\{\a_{n}^{*}\}_{n=0}^{\io}$
converging to $0$, define
$$ B^{*} = \sum_{n=0}^{\io} {1\over 2^{n}} \log {1 \over \a_{n}^{*}} ,
\qquad \G^{*}_{p} \= \sum_{n=0}^{\io} \a_{n}^{*} 2^{np} ,
\Eq(2.7) $$
for $p\in\ZZZ_{+}$. One has $\G^{*}_{p+1} > \G^{*}_{p}$
for all $p\in\ZZZ_{+}$ such that $\G^{*}_{p}$ is finite.}

\*

\0{\bf Lemma 2.}
{\it Let $\O\subset\RRR^{r}$ be an open set,
and let $\{\a_{n}^{*}\}_{n=0}^{\io}$ be a decreasing sequence
converging to zero such that one has $B^{*}<\io$
and $\G^{*}_{r}=C_{0}$ for some finite constant $C_{0}$.
Call $\O(C_{0})$ the subset of Bryuno vectors in $\BBBB_{r}(\O)$
such that $\a_{n}(\oo) \ge \a_{n}^{*}$ for all $n\ge 1$.
Then the Lebesgue measure of the set $\O^{\rm c}(C_{0})=
\O\setminus \O(C_{0})$ is bounded proportional to $C_{0}$.}

\*

\0{\it Proof.} The measure of the set $\O^{\rm c}(C_{0})$
can be bounded by
$$ \eqalign{
{\rm meas}(\O^{\rm c}(C_{0})) & \le
{\rm const.} \sum_{n=0}^{\io} \sum_{2^{n-1}<|\nn| \le 2^{n}}
{\a_{n}^{*} \over |\nn|} \cr
& \le {\rm const.} \sum_{n=0}^{\io}
\a_{n}^{*} 2^{nr} 2^{-(n-1)} \le {\rm const.} \, \G^{*}_{r-1}
\le {\rm const.} \, C_{0} , \cr}
\Eq(2.8) $$
so that the assertion follows. \qed

\*

The {\it Diophantine vectors}, that is the vectors
satisfying the usual Diophantine condition
$$ \left| \oo\cdot\nn \right| > { C_{0} \over |\nn|^{\t}} ,
\Eq(2.9) $$
for all $\nn\in\ZZZ^{r}_{*}$ and for suitable positive
constants $C_{0}$ and $\t$, are a particular case of Bryuno vectors,
with $\a_{n}(\oo) \ge 2^{-n\t} C_{0}$. In such a case in order to have
the convergence of the sum in \equ(2.8), hence to apply Lemma 2,
one must have $\t>r-1$, which is the condition
for the set of Diophantine vectors to have full measure.

The condition $\G^{*}_{r}=C_{0}$ motivates us to introduce
a new sequence $\{\g_{n}(\oo)\}_{n=1}^{\io}$,
with $\g_{n}(\oo)=C_{0}^{-1}\a_{n}(\oo)$, such that,
by setting $\a_{n}^{*}=C_{0}\g_{n}^{*}$ and defining
$$ \eqalign{
& \lis \G_{p}(\oo) = \sum_{n=0}^{\io} \g_{n}(\oo) 2^{np} , \qquad
\lis \G^{*}_{p} = \sum_{n=0}^{\io} \g_{n}^{*} 2^{np} = 1 , \cr
& \lis B(\oo) = \sum_{n=0}^{\io}
{1 \over 2^{n}} \log {1 \over \g_{n}(\oo) } , \qquad
\lis B^{*} = \sum_{n=0}^{\io}
{1 \over 2^{n}} \log {1 \over \g_{n}^{*}} , \cr}
\Eq(2.10) $$
one has $ \lis\G_{r}(\oo) \ge \lis\G^{*}_{r}$ and $\lis B(\oo) \le
\lis B^{*}$ for all $\O\subset \RRR^{r}$ and all $\oo \in \O(C_{0})$.

Note that if $|\oo\cdot\nn| < C_{0} \g_{n}(\oo)$ then $|\nn|>2^{n}$.
This is easily checked by contradiction: if $|\nn|\le 2^{n}$ then
$|\oo\cdot\nn| \ge C_{0} \g_{n}(\oo)$.

\*\*
\section(3,Fixing the perturbation parameter: proof of Theorem 1)

\0We follow very closely Ref.~\cita{GG2} (and Ref.~\cita{Ge3}), by
confining ourselves to show where the analysis differs. Also notations
which are not defined below are meant the same as in Ref.~\cita{GG2}.

The multiscale decomposition is performed as in Ref.~\cita{GG2},
by using the C$^{\io}$ non-decreasing function defined as
$$ \chi(x) = \cases{
1 , & if $|x|<C_{0}^{2}/4$, \cr
0 , & if $|x|>C_{0}^{2}$, \cr}
\Eq(3.1) $$
with the only difference that now $\chi_{n}$ for $n\ge 0$
is defined as $\chi_{n}(x) = \chi(\b^{-2}(\g^{*}_{n})^{-2}(\oo)x)$,
with $\b=1/4$. We set also $\psi_{n}(x)=1-\chi_{n}(x)$ for $n\ge 0$.

We define {\it clusters} and {\it self-energy clusters} as in
Ref.~\cita{GG2}, and we introduce the {\it self-energy value}
$\VV_{T}(x;\e,\oo)$ and the {\it tree value} $\Val(\th)$ according to
Ref.~\cita{GG2}, equations (5.8) and (5.11).
Of course one has ${\rm d}\VV_{T}/{\rm d}\oo =
\dpr_{x} \VV_{T} \dpr_{\oo} x + \dpr_{\oo}\VV_{T}$, and
${\rm d}\VV_{T}/{\rm d}\e = \dpr_{\e}\VV_{T}$.
Note that here and henceforth, with respect to Ref.~\cita{GG2},
we are making explicit the dependence of all quantities on $\oo$,
as we are interested also in changing $\oo$ for fixed $\e$.

In terms of the self-energy values we can define the
{\it self-energy matrices}
$$ \eqalign{
\MM^{[\le n]}(x;\e,\oo) & = \sum_{p=0}^{n} \MM^{[n]}(x;\e,\oo) , \cr
\MM^{[n]}(x;\e,\oo) & =
\left( \prod_{p=0}^{n} \chi_{p}(\D^{[p]}(x;\e,\oo)) \right) 
\sum_{T \in \SS^{\RR}_{k,n-1}} \VV_{T}(x;\e,\oo) , \cr}
\Eq(3.2) $$
where $\SS^{\RR}_{k,n}$ denotes the set of {\it renormalised
self-energy clusters} of degree $k$ and scale $[n]$ (renormalised
means that they do not contain any other self-energy clusters).
Such matrices are formally Hermitian (cf. Lemma 2 in Ref.~\cita{GG2}),
so that they admit $d$ real eigenvalues, which we denote by
$\l^{[n]}_{1}(x;\e,\oo)$, $\ldots$, $\l^{[n]}_{d}(x;\e,\oo)$.

The matrices $\MM^{[n]}(x;\e,\oo)$ in \equ(3.2) can be written as
$$ \MM^{[n]}(x;\e,\oo) = \left( \matrix{
\MM^{[n]}_{\a\a}(x;\e,\oo) & \MM^{[n]}_{\a\b}(x;\e,\oo) \cr
\MM^{[n]}_{\b\a}(x;\e,\oo) & \MM^{[n]}_{\b\b}(x;\e,\oo) \cr} \right) ,
\Eq(3.3) $$
where the labels $\a$ and $\b$ run over $\{1,\ldots,r\}$
and $\{r+1,\ldots,d\}$, respectively.

With respect to Ref.~\cita{GG2}, we slightly change the definition of
the {\it propagator divisors} for $n \ge 0$ (cf. Definition 6
in Ref.~\cita{GG2}); see also Ref.~\cita{GGG2}. We set
$$ \D^{[n]}(x;\e,\oo) = \left( {1 \over d} \sum_{j=1}^{d}
{1 \over ( x^{2}-\ul{\l}^{[n]}_{j}(\e,\oo) )^{2}}  \right)^{-1/2} ,
\Eq(3.4) $$
and define the {\it propagators} as
$$ g^{[n]}(x;\e,\oo) = \left( \prod_{p=0}^{n-1}
\chi_{p}(\D^{[p]}(x;\e,\oo)) \right) 
\psi_{n}(\D^{[n]}(x;\e,\oo))
\left( x^{2} - \MM^{[\le n]}(x;\e,\oo) \right)^{-1} ,
\Eq(3.5) $$
where the {\it self-energies} $\ul{\l}^{[n]}_{j}(\e,\oo)$ are
defined recursively as $\ul{\l}^{[n]}_{j}(\e,\oo)=\l^{[n]}_{j}
(\sqrt{\ul{\l}^{[n-1]}_{j}(\e,\oo)};\e,$ $\oo)$ for $n \ge 1$,
with $\ul{\l}^{[0]}_{j}(\e,\oo)=\e a_{j-r}$ for $j=r+1,\ldots,d$
and $\ul{\l}^{[0]}_{j}(\e,\oo)=0$ for $j=1,\ldots,r$.

Therefore if a line $\ell$ is on scale $[n]$, with $n\ge1$,
such that $g^{[n]}(\oo\cdot\nn_{\ell};\e,\oo) \neq 0$ one has
$$ \eqalign{
\min_{1\le j \le d} \left| (\oo\cdot\nn_{\ell})^{2} -
\ul{\l}^{[n]}_{j}(\e,\oo) \right|
& \ge {C_{0}^{2} \over 4\sqrt{d}} \b^{2} (\g^{*}_{n})^{2} , \cr
\min_{1\le j \le d} \left| (\oo\cdot\nn_{\ell})^{2} -
\ul{\l}^{[p]}_{j}(\e,\oo) \right|
& \le C_{0}^{2} \b^{2} (\g^{*}_{p})^{2} , \qquad 0 \le p \le n-1 , \cr}
\Eq(3.6) $$
and $\b$ is chosen in such a way to make uninfluential
the small changes of the propagators when shifting the lines
in order to exploit the cancellations discussed in Ref.~\cita{GG2},
Appendix A3 (see the proof of Lemma 5 below).
If a line $\ell$ is on scale $[0]$ and
$g^{[0]}(\oo\cdot\nn_{\ell};\e,\oo) \neq 0$ the condition
\equ(3.6) has to be replaced with $\min_{1\le j \le d}
|(\oo\cdot\nn_{\ell})^{2}-\ul{\l}^{[0]}_{j}(\e,\oo)|
\ge C_{0}^{2} \b^{2}(\g^{*}_{0})^{2}/4\sqrt{d}$.

We define the renormalised expansion
for $\hh(\pps,\bb_{0},\oo,\e)=(\aaa(\pps,\bb_{0},\oo,\e),
\bbb(\pps,\bb_{0},\oo,\e))$, by setting
$$ \eqalign{
& \hh(\pps,\bb_{0},\oo,\e) =
\sum_{\nn\in\zzz^{r}} {\rm e}^{i\nn\cdot\pps} \,
\hh_{\nn}(\bb_{0},\oo,\e) , \qquad
\hh_{\nn}(\bb_{0},\oo,\e) \= \hh_{\nn} , \cr
& \hh_{\nn} = (h_{\nn,1},\ldots,h_{\nn,d}) , \qquad h_{\nn,\g} =
\sum_{k=1}^{\io} \sum_{\th \in\Th^{\RR}_{k,\nn,\g}} \Val(\th)  , \cr}
\Eq(3.7) $$
where the set of trees $\Th^{\RR}_{k,\nn,\g}$ is defined
as in Ref.~\cita{GG2}, Definition 5.

We shall impose the following Diophantine conditions:
$$ \eqalign{
& \left| \oo\cdot\nn \pm \sqrt{\ul{\l}^{[n]}_{i}(\e,\oo)}
\right| \ge C_{0} \g^{*}_{n(\nn)} , \cr
& \left| \oo\cdot\nn \pm \sqrt{\ul{\l}^{[n]}_{i}(\e,\oo)} \pm
\sqrt{\ul{\l}^{[n]}_{j}(\e,\oo)} \right|
\ge C_{0} \g^{*}_{n(\nn)} \cr}
\Eq(3.8) $$
for all $i,j=1,\ldots,d$, for all $\nn\in\ZZZ^{r}_{*}$
and for all $n\ge 0$. We shall refer to conditions \equ(3.8)
as to the {\it first Mel$'\!$nikov conditions} (first line)
and the {\it second Mel$'\!$nikov conditions} (second line).

\*

\0{\bf Definition 4}
{\it Given $C_{0}\in\RRR_{+}$ and an open set
$\O\subset \RRR^{r}$ call $\O_{*}(C_{0}) \subset \O$
the set of Bryuno vectors in $\O(C_{0})$ satisfying
all the conditions \equ(3.8).}

\*

Hence the vectors $\oo\in\O_{*}(C_{0})$
verify the condition $\g_{n}(\oo) \ge \g_{n}^{*}$ for all
$n \ge 0$, with $C_{0}\g_{n}^{*}=\a_{n}^{*}$ and
the sequence $\{\a_{n}^{*}\}_{n=0}^{\io}$ defined as in Lemma 2,
and the first and second Mel$'$nikov conditions \equ(3.8).

\*

\0{\bf Lemma 3.} {\it Call $N_{n}(\th)$ the set of lines
in $\L(\th)$ which are on scale $[n]$. One has
$$ N_{n}(\th) \le K \, 2^{-n} M(\th) , \qquad
M(\th) = \sum_{\vvvv \in V(\th)} |\nn_{\vvvv}| ,
\Eq(3.9) $$
for a suitable constant $K$. One can take $K=2$.}

\*

\0{\it Proof.} First of all note that one can have $N_{n}(\th)$
only if $M(\th)\ge 2^{n-1}$. Indeed if a line $\ell$
is on scale $[n]$ then there exists $i\in\{1,\ldots,d\}$ such that
$$ C_{0} \g^{*}_{n-1} > C_{0} \b \g^{*}_{n-1} \ge
\left| | \oo\cdot\nn_{\ell} | -
\sqrt{\ul{\l}^{[n-1]}_{i}(\e,\oo)} \right| >
C_{0} \g^{*}_{n(\nn_{\ell})} ,
\Eq(3.10) $$
so that $n(\nn_{\ell}) \ge  n$. Then one must have
$|\nn_{\ell}|>2^{n(\nn_{\ell})-1} \ge 2^{n-1}$, hence
$M(\th) \ge |\nn_{\ell}| > 2^{n-1}$, thence
$K2^{-n}M(\th) \ge 1$ if $K\ge 2$.

Then one proves the bound $N_{n}(\th) \le \max\{2^{-n} K M(\th)-1,0\}$
for all $n\ge 0$, by induction on the number of vertices of the tree.
The only case which requires a different discussion with respect to
Ref.~\cita{GG2}, Appendix A3, is the one in which the root line $\ell$
is on scale $[n]$ and exits a cluster on scale $[n_{T}]$,
which has only one entering line, say $\ell'$, on scale $[n']$,
with $n' \ge n$. In such a case one has
$n_{T}<n$ of course, and, for suitable $i$ and $j$,
$$ \eqalign{
& \left| | \oo\cdot\nn | - \sqrt{\ul{\l}^{[n-1]}_{i}(\e,\oo)} \right| \le
C_{0} \b \g^{*}_{n-1} , \cr
& \left| | \oo\cdot\nn' | - \sqrt{\ul{\l}^{[n-1]}_{j}(\e,\oo)} \right| \le
C_{0} \b \g^{*}_{n-1} , \cr}
\Eq(3.11) $$
where $\nn=\nn_{\ell}$ and $\nn'=\nn_{\ell'}$, so that,
for suitable $\h,\h'\in\{\pm1\}$,
$$ \left| \oo\cdot (\nn-\nn') + \h \sqrt{\ul{\l}^{[n-1]}_{i}(\e,\oo)} +
\h' \sqrt{\ul{\l}^{[n-1]}_{j}(\e,\oo)} \right| < C_{0} \g^{*}_{n-1} ,
\Eq(3.12) $$
which by the Diophantine conditions \equ(3.8) implies
$n(\nn-\nn') \ge n$, hence one finds $M(T) \ge |\nn-\nn'|>2^{n-1}$,
if $M(T) = \sum_{\vvvv \in V(T)} |\nn_{\vvvv}|$. Call $\th'$ the
tree having $\ell'$ as root line. Then by the inductive hypothesis
$N_{n}(\th)=1+N_{n}(\th')\le 1 + \max\{2^{-n} K M(\th')-1,0\}$.
If the maximum is $0$ the bound is trivially satisfied,
because in such a case $N_{n}(\th)=1$ and we have seen that
in order to have a line on scale $[n]$ one needs $M(\th)>2^{n-1}$.
Otherwise one has $N_{n}(\th) \le 1 + 2^{-n} K M(\th')-1 \le
2^{-n} K M(\th) - 1 + (1 - 2^{-n} K M(T))$, where
$2^{-n} K M(T) \ge 1$ by the inequality $|\nn-\nn'|>2^{n-1}$,
provided that one takes $K\ge 2$. \qed

\*

\0{\bf Lemma 4.} {\it Call $N_{n}(T)$ the set of lines in $\L(T)$
which are on scale $[n]$, for $n \le n_{T}$. One has
$$ M(T) = \sum_{\vvvv \in V(T)} |\nn_{\vvvv}| > 2^{n_{T}-1} , 
\qquad N_{n}(T) \le K \, 2^{-n} M(T) ,
\Eq(3.13) $$
with the same constant $K$ as in \equ(3.9).}

\*

\0{\it Proof.} The first bound in \equ(3.13) can be proved by
{\it reductio ad absurdum} as in Ref.~\cita{GG2}, while the proof
of the second one is based on the  same argument used for proving
Lemma 3 (cf. Ref.~\cita{GG2}, Appendix A3, for further details). \qed

\*

Another difference with respect to Ref.~\cita{GG2} relies
in discussing the change of scale of the lines when performing
the cancellations inside the families $\FF_{T}$, when looking for
bounds on the entries of the matrices $\MM^{[n]}(x;\e,\oo)$.

\*

\0{\bf Lemma 5.} {\it Assume that the propagators
$g^{[p]}(x;\e,\oo)$ can be uniformly bounded for all $0 \le p\le n-1$ as
$$ \left| g^{[p]}(x;\e,\oo) \right| \le K_{1} C_{0}^{-2}
(\g^{*}_{p})^{-K_{2}} ,
\Eq(3.14) $$
for some $p$-independent constants $K_{1}$ and $K_{2}$.
Assume also that $\e$ is small enough.
Then, with the notations \equ(3.3), one has
$$  \eqalign{
\left\| \MM^{[n]}_{\a\a} (x;\e,\oo) \right\| & \le
B {\rm e}^{-\k_{1} 2^{n}} \min\{ \e^{2},\e x^{2} \} , \cr
\left\| \MM^{[n]}_{\a\b} (x;\e,\oo) \right\| & \le
B {\rm e}^{-\k_{1} 2^{n}} \min\{ \e^{2},\e^{3/2} |x| \} , \cr
\left\| \MM^{[n]}_{\b\b} (x;\e,\oo) \right\| & \le
B {\rm e}^{-\k_{1} 2^{n}} \e^{2} , \cr}
\Eq(3.15) $$
for suitable $n$-independent constants $B$ and $\k_{1}$.}

\*

\0{\it Proof.} Again we only discuss the differences with respect to
Ref.~\cita{GG2}. First we show that no cancellation is needed for
self-energy clusters $T$ with $C_{0}\g^{*}_{n(M(T))} \le 4|\oo\cdot\nn|$,
if $\nn$ is the momentum flowing through the entering line of $T$.
Note that we can extract from the self-energy value
a factor ${\rm e}^{-\k_{0} M(T)/4} \le {\rm e}^{-\k_{0} 2^{n(M(T))}/8}$.
If we set $2^{-n}\log 1/\a_{n}^{*}=a_{n}$, we have $\lim_{n\to\io}a_{n}=0$
(because $B^{*}<\io$), hence for $\oo\cdot\nn$ small enough ${\rm e}^{-\k_{0}
2^{n(M(T))}/8} \le (C_{0}\g^{*}_{n(M(T))})^{\k_{0}/8a_{n(M(T))}}
\le C_{0}^{2} (\g^{*}_{n(M(T))})^{2} \le 16|\oo\cdot\nn|^{2}$. 

Hence we need the cancellations only for self-energy clusters $T$ with
$C_{0}\g^{*}_{n(M(T))} > 4|\oo\cdot\nn|$ if $\nn$ is the momentum
flowing through the entering line of $T$. In such a case
one can reason as follows. For any line $\ell\in \L(T)$
and for any $n \le n_{\ell}$ one has,
by the Diophantine conditions \equ(3.8),
$||\oo\cdot\nn_{\ell}^{0}|-\sqrt{\ul{\l}^{[n]}_{j}(\e,\oo)}|
\ge C_{0}\g^{*}_{n(\nn_{\ell}^{0})}$, if $\nn_{\ell}^{0}$ is defined as
$$ \nn_{\ell}^{0} = \sum_{\wwwww \in V(T) \atop \wwwww \le \vvvvv}
|\nn_{\vvvv}| , \qquad \ell \= \ell_{\vvvv} .
\Eq(3.16) $$
On the other hand one has $|\nn_{\ell}^{0}| \le M(T)$, so that
$C_{0}\g^{*}_{n(\nn_{\ell}^{0})} \ge C_{0}\g^{*}_{n(M(T))} >
4|\oo\cdot\nn|$. Therefore we can bound
$$ 2 \left| |\oo\cdot\nn_{\ell}^{0}| -
\sqrt{\ul{\l}^{[n]}_{j}(\e,\oo)} \right|
\ge \left| |\oo\cdot\nn_{\ell}| -
\sqrt{\ul{\l}^{[n]}_{j}(\e,\oo)} \right|
\ge {1\over 2}
\left| | \oo\cdot\nn_{\ell}^{0} | -
\sqrt{\ul{\l}^{[n]}_{j}(\e,\oo)} \right| .
\Eq(3.17) $$
This implies the following. When considering a family $\FF_{T}$, a line
$\ell\in \L(T)$ with momentum $\nn_{\ell}$ can be on a scale $[n_{\ell}]$
such that $g^{[n_{\ell}]}(\oo\cdot\nn_{\ell}^{0};\e,\oo) = 0$.
But in such a case it is obtained, by shifting the external
lines of $T$, from a line with non-vanishing propagator,
that is from a line for which \equ(3.6) holds.
Then, even if the bounds \equ(3.6) can fail to hold, one still
obtains bounds of the same form with the only difference
that $\b^{2}$ is replaced with $\b^{2}/4$ in the first line
and with $4\b^{2}$ in the second line. In particular for $\b=1/4$
the inequalities \equ(3.10) and \equ(3.12) are still satisfied
for all self-energy clusters in the family $\FF_{T}$.

This shows that one can reason in Ref.~\cita{GG2} to
deduce the bounds \equ(3.15), which are of algebraic nature,
and are due to symmetry properties of the self-energy matrices,
that is $\MM^{[\le n]T}(x;\e,\oo)=\MM^{[\le n]}(-x;\e,\oo)$
and $\MM^{[\le n]\dagger}(x;\e,\oo)=\MM^{[\le n]}(x;\e,\oo)$. \qed

\*

\0{\bf Lemma 6.} {\it Assume that the propagators
$g^{[p]}(x;\e,\oo)$ can be uniformly bounded for all $0 \le p\le n-1$
as in \equ(3.14), for some $p$-independent constants $K_{1}$ and $K_{2}$.
Assume also that $\e$ is small enough.
The matrices $\MM^{[\le n]}(x;\e,\oo)$ are differentiable in $x$ in
the sense of Whitney, and for all $x',x$ where they are defined one has
$$ \eqalignno{
& \left\| \MM^{[\le n]}(x';\e,\oo) -
\MM^{[\le n]}(x;\e,\oo) - \dpr_{x}
\MM^{[\le n]}(x;\e,\oo) \, (x'-x) \right\| =
\e^{2} \, o(|x'-x|) , \cr
& \left\| \dpr_{x} \MM^{[\le n]}(x;\e,\o) \right\|
\le B \e^{2} ,
& \eq(3.18) \cr} $$
for a suitable positive constant $B$. Moreover for all $j=1,\ldots,d$
and for a suitable constant $A$ one has
$$ \left| \dpr_{x} \l^{[\le n]}_{j}(x;\e,\oo) \right| \le A \e^{2} ,
\Eq(3.19) $$
where the derivative is in the sense of Whitney.}

\*

\0{\it Proof.} The proof of \equ(3.18) can be performed
as in Ref.~\cita{GG2}; cf. in particular Section 6.
Then property \equ(3.19) follows from general properties
of Hermitian matrices. One can refer to Ref.~\cita{GG2}, Appendix A4,
in the case in which the eigenvalues $a_{i}$ are all distinct.
Otherwise one can apply the results on non-analytic
Hermitian matrices discussed in Ref.~\cita{Ka}, Chapter 2, Section 6:
one can rely on Rellich's theorem \cita{Re} to deduce differentiability
of the eigenvalues and on Lidski\u\i's theorem \cita{L}
to obtain a bound on the derivative. \qed

\*

\0{\bf Lemma 7.} {\it Assume that the propagators
$g^{[p]}(x;\e,\oo)$ can be uniformly bounded for all $0 \le p\le n-1$
as in \equ(3.14), for some $p$-independent constants $K_{1}$ and $K_{2}$.
Assume also that $\e$ is small enough.
The self-energies $\ul{\l}^{[p]}_{j}(\e,\oo)$ satisfy
for all $0 \le p\le n$ and all $1\le j \le d$ the closeness property
$$ \left| \ul{\l}^{[p]}_{j}(\e,\oo) - \ul{\l}^{[p-1]}_{j}(\e,\oo) \right|
\le A {\rm e}^{-\k_{1} 2^{p}} \e^{2}
\Eq(3.20) $$
and one has
$$ \left| \l^{[p]}_{j}(x;\e,\oo) \right| \le
A \min\{\e^{2},\e x^{2} \} ,
\qquad j=1,\ldots,r ,
\Eq(3.21) $$
for a suitable positive constant $A$.}

\*

\0{\it Proof.} The proof can be performed as in Ref.~\cita{GG2}, by using
the bounds \equ(3.15) and the fact that the matrices $\MM^{[\le n]}
(x;\e,\oo)$ are Hermitian (see Ref.~\cita{GG2}, Lemma 2).

\*

\0{\bf Lemma 8.} {\it Assume that the propagators
$g^{[p]}(x;\e,\oo)$ can be uniformly bounded for all $0 \le p\le n-1$
as in \equ(3.14), for some $p$-independent constants $K_{1}$ and $K_{2}$.
Assume also that $\e$ is small enough.
If $g^{[n]}(x;\e,\oo) \neq 0$ then one has
$$ \min_{j=1,\ldots,d}
\left| x^{2} - \l^{[n]}_{j}(x;\e,\oo) \right| \ge
{1 \over 2} \min_{j=1,\ldots,d}
\left| x^{2} - \ul{\l}^{[n]}_{j}(\e,\oo) \right| .
\Eq(3.22) $$
The same holds if $g^{[n]}(x;\e,\oo) = 0$ but \equ(3.6)
are satisfied with $\b^{2}$ replaced with $\b^{2}/4$
in the first line and with $4\b^{2}$ in the second line.}

\*

\0{\it Proof.} The inequality \equ(3.22) can be proved
by induction on $n$. For $n=0$ it is trivially satisfied
as the matrix $\MM^{[0]}(x;\e,\oo)$ does not depend on $x$.
Let us assume that the inequality holds for all $0 \le p\le n-1$,
and let us show that then it follows also for $p=n$.
First of all note that in this case we can apply Lemma 7,
so that one has
$$ \left| \ul{\l}^{[n]}_{j}(\e,\oo) - \ul{\l}^{[n-1]}_{j}(\e,\oo) \right|
\le A {\rm e}^{-\k_{1} 2^{n}} \e^{2} \le A_{1} \e^{2} (\g^{*}_{n})^{2} ,
\Eq(3.23) $$
for some constant $A_{1}$ (we have used $\oo\in\BBBB_{r}$ to deduce
${\rm e}^{-\k_{1} 2^{n}} \le {\rm const.} C_{0}^{2} (\g^{*}_{n})^{2}$).
By hypothesis one has $g^{[n]}(x;\e,\oo) \neq 0$, hence, by \equ(3.6),
$$ \left| x^{2}-\ul{\l}_{j}^{[n]}(\e,\oo) \right| \ge
C_{1} (\g^{*}_{n})^{2} ,
\Eq(3.24) $$
for all $1 \le j \le d$ and some constant $C_{1}$.
One the other hand we can write for all $j'=1,\ldots,d$
$$ \left| x^{2} - \l^{[n]}_{j'}(x;\e,\oo) \right|
\ge \left| x^{2} - \ul{\l}^{[n]}_{j'}(\e,\oo) \right| -
\left| \l^{[n]}_{j'}(x;\e,\oo) - \ul{\l}^{[n]}_{j'}(\e,\oo) \right| ,
\Eq(3.25) $$
where we recall that, by construction, $\ul{\l}^{[n]}_{j'}(\e,\oo)=
\l^{[n]}_{j'}(\sqrt{\ul{\l}^{[n-1]}_{j'}(\e,\oo)};\e,\oo)$. Hence if
$j'>r$ and $x>0$ we can bound, by Lemma 6 and in particular \equ(3.19),
$$ \eqalign{
\left| \l^{[n]}_{j'}(x;\e,\oo) - \ul{\l}^{[n]}_{j'}(\e,\oo) \right| & \le
A \e^{2} \left| |x| - \sqrt{\ul{\l}^{[n-1]}_{j'}(\e,\oo)} \right| \cr
& \le A' \e^{3/2} \left| x^{2} - \ul{\l}^{[n-1]}_{j'}(\e,\oo) \right|
\le 2A' \e^{3/2} \left| x^{2} - \ul{\l}^{[n]}_{j'}(\e,\oo) \right| , \cr}
\Eq(3.26) $$
where we have used the relations \equ(3.23) and \equ(3.24).
If $j'>r$ and $x<0$ we can apply the same argument by using
the symmetry property that $\l^{[n]}_{j'}(-x;\e,\oo)$
belongs to the spectrum if $\l^{[n]}_{j'}(x;\e,\o)$ does
(because $\MM^{[\le n]}(-x;\e,\oo)=(\MM^{[\le n]}(x;\e,\oo))^{T}$,
see Lemma 2, (ii), in Ref.~\cita{GG2}; cf. also the comments
at the end of the proof of Lemma 5 above). If $j' \le r$ then
$$ \left| \l^{[n]}_{j'}(x;\e,\oo) - \ul{\l}^{[n]}_{j'}(\e,\oo) \right| =
\left| \l^{[n]}_{j'}(x;\e,\oo) \right| \le A \e x^{2} =
A \e \left( x^{2} - \ul{\l}^{[n]}_{j'}(\e,\oo) \right) ,
\Eq(3.27) $$
by \equ(3.21). By inserting \equ(3.26) or \equ(3.27) into \equ(3.25),
and choosing $j'$ as the value of $j$ minimising $|x^{2}-\l^{[p]}_{j}
(x;\e,\oo)|$, then the assertion follows. \qed

\*

\0{\bf Lemma 9.} {\it Let $\oo\in\O_{*}(C_{0})$ and assume that $\e$
is small enough, say $|\e|<\e_{0}$. Then the series
\equ(3.7) admits the bound $|h_{\nn,\g}| < H \,
{\rm e}^{-\k |\nn|} |\e|$ for suitable positive constants $\k$ and $H$.
One has $\e_{0}=O(C_{0}^{2}(\g^{*}_{m_{0}})^{2})$ with $m_{0}$
depending on $\k_{0}$.}

\*

\0{\it Proof.} One can proceed as in Ref.~\cita{GG2}.
Here we outline only the differences. As a consequence of Lemma 8,
we can prove inductively that for all $n\ge0$ the propagators
with scales $[n]$ are bounded proportionally to $(\g_{n}^{*})^{-K_{2}}$,
and one finds, in particular, $K_{2}=2$. Then
the product of propagators can be bounded by relying on Lemma 3
for the lines on scale $[n]$, with $n> m_{0}$,
and bounding with $(\g^{*}_{m_{0}})^{-2k}$
the propagators of all lines on scale $[n]$ for $n \le m_{0}$.
Therefore for any tree $\th\in\Th^{\RR}_{k,\nn,\g}$ one has
$$ \eqalign{
\prod_{\ell \in \L(\th)} \left| g^{[n_{\ell}]} \right|
& \le \left( {\rm const.} \right)^{k}
C_{0}^{-2k} \left( {1 \over \g^{*}_{m_{0}} } \right)^{2k}
\prod_{n=m_{0}+1}^{\io}
\left( {1 \over \g^{*}_{n}} \right)^{2N_{n}(\th)} \cr
& \le \left( {\rm const.} \right)^{k}
C_{0}^{-2k} \left( {1 \over \g^{*}_{m_{0}} } \right)^{2k}
\exp \left( K |\nn| \sum_{n=m_{0}+1}
{1 \over 2^{n}} \log {1 \over \g^{*}_{n}} \right) , \cr}
\Eq(3.28) $$
and one can choose $m_{0}=m_{0}(\k_{0})$,
so the last exponential is less than ${\rm e}^{\k_{0}|\nn|/4}$.
By making use of the bound
$\prod_{\vvvv\in V(\th)} {\rm e}^{-\k|\nn_{\vvvv}|} \le
{\rm e}^{-\k|\nn|}$, this produces an overall
factor ${\rm e}^{-\k_{0}|\nn|/2}$. This completes the proof,
and it gives $\e_{0}=O(C_{0}^{2}(\g^{*}_{m_{0}})^{2})$. \qed

\*

Note that for Diophantine vectors satisfying the bound \equ(2.8) one has
$m_{0}=\t \, O(\log 1/\k_{0})$, and one obtains $\e_{0}=O(C_{0}^{2})$,
for fixed $\t$ and $C_{0}$.

If we are interested in studying the conservation of a maximal torus
with rotation vector $\oo$ satisfying the Bryuno condition
$B(\oo)<\io$, we can use directly the sequence
$\{\g_{n}(\oo)\}_{n=0}^{\io}$ for the multiscale decomposition,
without introducing a further sequence $\{\g_{n}^{*}\}_{n=0}^{\io}$.
Then the result stated in Lemma 9 holds with $\g_{m_{0}}(\oo)$
replacing $\g_{m_{0}}^{*}$.

An important remark is that, in the case of perturbations
which are trigonometric polynomials of degree $N$ in the bound
\equ(3.9) one can bound $M(\th) \le k N$, and as
consequence the product of propagators in \equ(3.28) can be bounded as
$$ \prod_{\ell\in\L(\th)} \left| g^{[n_{\ell}]}_{\ell} \right|
\le \exp \left( 2 K N k \sum_{n=0}^{\io} {1 \over 2^{n}} \log
{1 \over \a_{n}(\oo)} \right) = {\rm e}^{4N k B(\oo)} ,
\Eq(3.29) $$
which implies $\e_{0}=O({\rm e}^{-4N B(\oo)})$. We can compare this
result with the one found in Ref.~\cita{GL} for maximal tori,
where a bound of this kind with the factor $4$ replaced with
the likely optimal $2$ was obtained. With the techniques described
in this paper some further work is necessary in order to reach
the factor $2$; cf. for instance Ref.~\cita{BeG}.
On the other hand an advantage with respect
to Ref.~\cita{GL}, which relies on using Lie transforms for
Hamiltonian flows, is that our techniques apply, essentially unchanged,
not only to the case of flows, but also to the case of diffeomorphisms,
as done for instance in Refs.~\cita{BeG} and \cita{Ge1},
where the case of the standard map was explicitly treated.

\*

\0{\bf Lemma 10.} {\it For all $\oo,\oo'\in \O_{*}(C_{0})$ one has
$$ \eqalign{
& \left\| \MM^{[\le n]}(\oo'\cdot\nn;\e,\oo') -
\MM^{[\le n]}(\oo\cdot\nn;\e,\oo) \right. - \cr
& \qquad \qquad \left. \dpr_{\oo}
\MM^{[\le n]}(\oo\cdot\nn;\e,\oo) \cdot (\oo'-\oo) \right\| =
\e^{2} |\nn| \, o(|\oo'-\oo|) , \cr
& \left\| \dpr_{\oo} \MM^{[\le n]}(x;\e,\o) \right\|
\le B \e^{2}|\nn| , \cr}
\Eq(3.30) $$
for a suitable constant $B$, and, as a consequence,
$$ \left| \dpr_{\oo} \l^{[\le n]}_{j}(x;\e,\o) \right|
\le A \e^{2}|\nn| , \qquad
\left| \dpr_{\oo} \ul{\l}^{[\le n]}_{j}(\e,\o) \right|
\le A \e^{2}|\nn| ,
\Eq(3.31) $$
for a suitable constant $A$.}

\*

The proof is deferred to Appendix \secc(A1), which represents
the core of the technical part of the paper (by assuming
the results of Ref.~\cita{GG2} for granted). Note that in fact
it is enough to prove \equ(3.30), because then property \equ(3.31)
follows by general properties of Hermitian matrices;
cf. the comments in the proof of Lemma 6.

\*

\0{\bf Lemma 11.} {\it The Lebesgue measure of the set
$\O\setminus\O_{*}(C_{0})$ is bounded proportionally to $C_{0}$.}

\*

\0{\it Proof.} Let us start with the first conditions in \equ(3.8).
We can reason as in Ref.~\cita{Ge3}, and write $\oo=\a \nn/|\nn| + \bb$,
with $\bb\cdot\nn=0$. Then we define $\a(t)$, $t\in[-1,1]$, such that
$ F(t) = \a(t) |\nn| \pm \sqrt{\ul{\l}^{[n]}_{i}(\e,(\a(t),\bb))} =
t C_{0} \g_{n(\nn)}^{*}$, so that ${\rm d} F/{\rm d} t=|\nn|(1+O(\sqrt{\e}))
{\rm d}\a/{\rm d}t=C_{0}\g_{n(\nn)}^{*}$; cf. the proofs of Lemma 3 and
Lemma 4 in Ref.~\cita{Ge3} for further details. Given $p\ge 1$
we define $\O^{[p]}$ as the sets of $\oo\in\O$ satisfying the
conditions \equ(3.8) for $n \le p$; we also set $\O^{[0]}=\O$.
For each $n$, for fixed $\nn$ and $i$, we call $I_{n}(i,\nn)$
the sets of $\oo\in\O^{[n-1]}$ such that either
$| \oo\cdot\nn \pm \sqrt{\ul{\l}^{[n]}_{i}(\e,\oo)}| < C_{0}
\g^{*}_{n(\nn)}$ or $| \oo\cdot\nn \pm \sqrt{\ul{\l}^{[n]}_{i}(\e,\oo)}
\pm \sqrt{\ul{\l}^{[n]}_{j}(\e,\oo)}| < C_{0} \g^{*}_{n(\nn)}$.
We define in the same way the sets $J_{n}(i,\nn)$, with the only
difference that the width of the sets is $2C_{0} \g^{*}_{n(\nn)}$
instead of $C_{0} \g^{*}_{n(\nn)}$. By the closeness property of Lemma 7,
there is some $n_{1}(\nn) = O(\log\log 1/\g_{n(\nn)}^{*})$ such that
all the sets $I_{n}(i,\nn)$ fall inside $J_{n_{1}(\nn)}$ for
$n\ge n_{1}(\nn)$. Therefore for all $\nn\in\ZZZ_{*}^{r}$, all
$i=r+1,\ldots,d$, and all $n\le n_{1}(\nn)$ we have to exclude all
values of $\oo\in\O^{[n-1]}$ which fall inside the set $J_{n}(i,\nn)$;
we refer to Ref.~\cita{Ge2}, Section 7, for details.
Note that $\oo\in\BBBB_{r}$ implies $n_{1}(\nn) \le C n(\nn)$,
for some constant $C$. Hence we can bound the measure of the set
of excluded values by a constant times
$$ \eqalign{
{\rm const.} \sum_{i=r+1}^{d} C_{0} \sum_{\nn\in\zzz^{r}}
\sum_{n=1}^{n_{1}(\nn)} {\g^{*}_{n(\nn)} \over |\nn|}
& \le {\rm const.} \, s C_{0}
\sum_{n=0}^{\io} 2^{n(r-1)} \g_{n}^{*} \log\log 1/\g_{n}^{*} \cr
& \le {\rm const.} \,
s C_{0} \sum_{n=0}^{\io} n 2^{n(r-1)} \g_{n}^{*} , \cr}
\Eq(3.32) $$
which is bounded proportionally to $C_{0}$ by Lemma 2.

Analogously one discusses the other conditions in \equ(3.8). Simply one
defines $F(t) = \a(t) |\nn| \pm \sqrt{\ul{\l}^{[n]}_{i}(\e,(\a(t),\bb))}
\pm \sqrt{\ul{\l}^{[n]}_{j}(\e,(\a(t),\bb))} = t C_{0} \g_{n(\nn)}^{*}$,
so that again one has ${\rm d} F/{\rm d} t =|\nn|(1+O(\sqrt{\e}))
{\rm d}\a/{\rm d}t$ $=$ $C_{0}\g_{n(\nn)}^{*}$,
and one can proceed as before. \qed

\*

To complete the proof of Theorem 1 we have to prove the last
assertion about maximal tori, that is that for $r=d$
most of phase space is filled by invariant tori.

We summarise what we have found so far. The invariant tori are determined
by the corresponding rotation vectors $\oo$. For $\oo\in\O_{*}(C_{0})$
we can parameterise the invariant torus with rotation vector $\oo$ as
$$ \cases{
\aa = \pps + \aaa (\pps,\oo,\e) , \cr
\AA = \oo + (\oo\cdot\dpr_{\pps})\, \aaa (\pps,\oo,\e) , \cr}
\Eq(3.33) $$
with $\pps\in \TTT^{r}$. Moreover the function $\aaa$ is
analytic in $\pps$ (as the Fourier coefficients decay exponentially),
while it is defined only on a Cantorian set of values $\oo$.

For each value of $\pps$ we can consider the map $\oo \mapsto \AA(\oo)$,
defined in \equ(3.33). For $\e$ small enough in \equ(3.33) one has
$|\AA-\oo| = |(\oo\cdot\dpr_{\pps})\,\aaa(\pps,\oo,\e)| \le R|\e|$
for some constant $R$. Call $\O_{*}(C_{0},d)$ the open set obtained
from $\O_{*}(C_{0})$ by excluding all vectors within a distance $d$
from the boundary of $\O$, i.e. $\O_{*}(C_{0},d)=\{\oo\in\O(C_{0}):
d(\oo,\dpr\O) \ge d\}$, with $d(\oo,\dpr\O)=\min_{\oo'\in\dpr\O}
|\oo-\oo'|$, and define $\AAA_{*}(C_{0})$ as the image of
$\O_{*}(C_{0},R|\e|)$ of such a map (note that the latter is not just
the inverse of the frequency map, rather it is a perturbation of it).
Then the measure of the complement of the
action variable space filled by invariant tori is given by
$$ {\rm meas}(\AAA_{*}^{c}(C_{0})) = \int_{\AAA^{c}_{*}(C_{0})}
{\rm d}\AA = \int_{\O_{*}^{c}(C_{0},R|\e|)} {\rm d}\oo
\left| \det \dpr_{\oo}\AA \right| ,
\Eq(3.34) $$
provided that the Jacobian in the last integral is well defined
(that is the map $\oo\to \AA(\oo)$ is smooth enough,
at least in the sense of Whitney) and is uniformly bounded.
This turns out to be the case, as the following result shows.

\*

\0{\bf Lemma 12.} {\it The solutions of the equations of motion
$\hh(\pps,\bb_{0},\oo,\e)$ are differentiable in $\oo$
in the sense of Whitney for $\oo\in\O_{*}(C_{0},R|\e|)$.}

\*

\0{\it Proof.} The proof (for any value of $r \le d$,
not necessarily $r=d$) can be performed as
for Lemma 10, with the only difference that now we have to deal with
the renormalised expansion for $h_{\nn,\g}$ instead
of the matrices $\MM^{[\le n]}(x;\e,\oo)$, hence
with trees instead of self-energy clusters. The condition
$d(\oo,\dpr\O)\ge R|\e|$ yields that the actions variables $\AA$
remain in $\O$ for all values of $\pps$. \qed

\*

As a consequence we can bound the Jacobian in \equ(3.34)
by using \equ(3.33), which gives $\dpr_{\oo}\AA = \openone
+ \dpr_{\oo} (\oo\cdot\dpr_{\pps})\, \aaa (\pps,\oo,\e)$, and Lemma 12,
which assures that the last derivative (in the sense of Whitney)
is bounded proportionally to $\e$. Therefore we can bound
${\rm meas}(\AAA_{*}^{c}(C_{0}))$ proportionally to $C_{0}$
by Lemma 11, and by taking $C_{0}=O(\sqrt{|\e|})$ (which is allowed
by Lemma 9), we obtain the last assertion of Theorem 1.
Cf. also Refs.~\cita{CG} and \cita{P1}, where the usual
Diophantine conditions were considered
in the analytic and differentiable case, respectively.

\*\*
\section(4,Fixing the rotation vector: proof of Theorem 2)

\0In the following we assume that $\oo$ is fixed, and that it satisfies
the Bryuno Diophantine condition $B(\oo)<\io$, with
$B(\oo)$ defined in \equ(2.4). Set $\g_{n}=C_{0}^{-1}\a_{n}(\oo)$.

Let $\e\in(\e_{0}/4,\e_{0}]$ and set $\l^{[0]}_{d}=\e a_{s}$
and $\e_{0} a_{s} =\L_{0}$.
Define $n_{0}\in \NNN$ such that $C_{0}\g_{n_{0}+1} < 2 \sqrt{\L_{0}}
\le C_{0}\g_{n_{0}}$. We set
$$ \g_{n}^{*} = \cases{
\g_{n} , & $n < n_{0}$, \cr
\g_{n}2^{-n(r+1)} , & $n \ge n_{0}$, \cr}
\Eq(4.1) $$
and, by using the sequence $\{\g_{n}^{*}\}_{n=0}^{\io}$,
we proceed as in Section \secc(3), for constructing the
multiscale decomposition of the propagators. Though,
we define $\D^{[n]}(x;\e,\oo)=\D^{[0]}(x;\e,\oo)$ for $n \le n_{0}$.

The main difference is that 
we shall need the following Diophantine conditions:
$$ \eqalign{
& \left| \oo\cdot\nn \pm \sqrt{\ul{\l}^{[n]}_{i}(\e,\oo)}
\right| \ge C_{0} \g^{*}_{n(\nn)} , \cr
& \left| \oo\cdot\nn \pm \sqrt{\ul{\l}^{[n]}_{i}(\e,\oo)} \pm
\sqrt{\ul{\l}^{[n]}_{j}(\e,\oo)} \right| \ge
C_{0} \g^{*}_{n(\nn)} \cr}
\Eq(4.2) $$
for all $i,j=1,\ldots,d$, for all $\nn\in\ZZZ^{r}_{*}$
such that $n(\nn) \ge n_{0}$ and for all $n\ge n_{0}$.
We do not impose any conditions like \equ(4.2) for $n \le n_{0}$,
because for such scales one has $|\oo\cdot\nn| > 2\sqrt{\L_{0}}$,
so that we can bound $|(\oo\cdot\nn)^{2} - \ul{\l}^{[n]}_{i}(\e,\oo)|$
with $(\oo\cdot\nn)^{2}/2$ for all $i=1,\ldots,d$. In the same way
we have excluded in \equ(4.2) the values of $\nn\in\ZZZ^{r}_{*}$
such that $n(\nn) \le n_{0}$. Hence, at the price of adding a
factor $2^{2}$ in the bound of each propagator, we can confine
ourselves to impose \equ(4.3) only for $\nn$ such that
$n(\nn) \ge n_{0}$ and for $n\ge n_{0}$.

Then we can prove the following result.

\*

\0{\bf Lemma 13.} {\it Call $N_{n}(\th)$ the set of lines
in $\L(\th)$ which are on scale $[n]$. One has
$$ N_{n}(\th) \le K \, 2^{-n} M(\th) , \qquad
M(\th) = \sum_{\vvvv \in V(\th)} |\nn_{\vvvv}| ,
\Eq(4.3) $$
for a suitable constant $K$. One can take $K=2$.}

\*

\0{\it Proof.} The proof proceeds exactly as for Lemma 3,
with the only difference that we have to deal in a different way
with the lines on scales $n<n_{0}$ and those with
scales $n \ge n_{0}$. The same was done in Ref.~\cita{GG2}.\qed

\*

In the same way the following result is proved.

\*

\0{\bf Lemma 14.} {\it Call $N_{n}(T)$ the set of lines
in $\L(T)$ which are on scale $[n]$. One has
$$ M(T) = \sum_{\vvvv \in V(T)} |\nn_{\vvvv}| > 2^{n_{T}-1} , \qquad
N_{n}(T) \le K \, 2^{-n/2} M(T) ,
\Eq(4.4) $$
for a suitable constant $K$.}

\*

Then Lemma 5 is replaced with the following one.

\*

\0{\bf Lemma 15.} {\it
Assume that the propagators
$g^{[p]}(x;\e,\oo)$ can be uniformly bounded for
all $0 \le p\le n-1$ as \equ(3.14),
for some $p$-independent constant $K$.
Assume also that $\e_{0}$ is small enough.
With the notations \equ(3.3) one has
$$  \eqalign{
\left\| \MM^{[n]}_{\a\a} (x;\e,\oo) \right\| & \le
B {\rm e}^{-\k_{1} 2^{n/2}} \min\{ \e^{2},\e x^{2} \} , \cr
\left\| \MM^{[n]}_{\a\b} (x;\e,\oo) \right\| & \le
B {\rm e}^{-\k_{1} 2^{n/2}} \min\{ \e^{2},\e^{3/2} |x| \} , \cr
\left\| \MM^{[n]}_{\b\b} (x;\e,\oo) \right\| & \le
B {\rm e}^{-\k_{1} 2^{n/2}} \e^{2} , \cr}
\Eq(4.5) $$
for all $n\in\NNN$ and for suitable $n$-independent constants
$B$ and $\k_{1}$.}

\*

Therefore we can prove the following estimates. The proof
is the same as for Lemma 6, as one easily realizes
that it works for fixed values of $\e$ and $\oo$.

\*

\0{\bf Lemma 16.} {\it Let $\oo$ satisfy the Diophantine condition
\equ(1.4) and assume that $\e$ is small enough, say $|\e|<\e_{0}$.
Then the series \equ(3.7) admits the bound $|h_{\nn,\g}| < A \,
{\rm e}^{-\k |\nn|} |\e|$ for suitable positive constants $\k$ and $A$.
One has $\e_{0}=O(C_{0}^{2}(\g^{*}_{m_{0}})^{2})$ with $m_{0}$
depending on $\oo$ and $\k_{0}$.}

\*

With respect to Section 3 
the first differences appear when dealing with Whitney extensions
of the matrices $\MM^{[\le n]}(x;$ $\e,\oo)$: indeed now $\oo$
is assumed to be fixed, while $\e$ is the free parameter.
We define $\EE_{n_{0}}\=(\e_{0}/4,\e_{0}]$ and for $n >n_{0}$,
recursively, $\EE_{n}=\EE_{n-1}\setminus\EE_{n}^{o}$, where
$\EE_{n}^{o}$ is the set of values of $\e\in\EE_{n}$ such that
the conditions \equ(4.2) are violated. We define also
$\EE_{*} = \cap_{n=n_{0}}^{\io} \EE_{n}$.

\*

\0{\bf Lemma 17.} {\it For all $n\ge 0$ and all $\e,\e'\in \EE_{n}$
one has
$$ \eqalign{
& \left\| \MM^{[\le n]}(x;\e',\oo) - \MM^{[\le n]}(x;\e,\oo) - \dpr_{\e}
\MM^{[\le n]}(x;\e,\oo) \, (\e'-\e) \right\| = \e \, o(|\e'-\e|) , \cr
& \left\| \dpr_{\e} \MM^{[\le n]}(x;\e,\oo) \right\|
\le B , \cr}
\Eq(4.6) $$
and, as a consequence,
$$ \eqalign{
& B ' \le \left| \dpr_{\e} \ul{\l}_{j}^{[n]}(\e,\oo) \right| \le B ,
\qquad r+1 \le j \le d , \cr
&  B ' \le \left| \dpr_{\e} \left(
\ul{\l}_{i}^{[n]}(\e,\oo) \pm \ul{\l}_{j}^{[n]}(\e,\oo) \right)
\right| \le B ,
\qquad r+1 \le j<i \le d , \cr}
\Eq(4.7) $$
for suitable positive constants $B$ and $B'$.}

\*

\0{\it Proof.} The proof can be performed as for Lemma 10,
with the parameter $\e$ now playing the role of the
parameters $\oo$. We do not give the details, which,
however, have been worked out in Ref.~\cita{GG2}.
Again the upper bound \equ(4.7) follows from \equ(4.6);
cf. analogous comments in the proof of Lemma 6.
To obtain the lower bound one has to use also that
$\l_{j}(x;\e,\oo) = a_{j}\e + O(\e^{2})$,
with $a_{j} \neq 0$ and $a_{i} \neq a_{j}$ for
$i,j=r+1,\ldots,d$; again cf. Ref.\cita{GG2} for details. \qed

\*

Lemma 17 implies that the matrices
$\MM^{[\le n]}(x;\e,\oo)$ can be extended in $(0,\e_{0})$
to smooth $C^{1}$ functions (Whitney extensions). Again a closeness
property of the self-energies, which reads
$$ \left| \ul{\l}^{[n]}_{j}(\e,\oo) -
\ul{\l}^{[n-1]}_{j}(\e,\oo) \right|
\le B {\rm e}^{-\k_{1} 2^{n/2}} \e^{2} ,
\Eq(4.8) $$
follows from Lemma 15. As before the bounds \equ(4.8) can be
improved for the first $r$ self-energies, and give
$|\l^{[n]}_{j}(x;\e,\oo) | \le A \min\{ \e^{2}, \e x^{2} \}$,
$j=1,\ldots r$. What really changes with respect to the previous case
is the estimate of the set of allowed values of $\e$, which explains
why we have required the stronger condition on $\dpr_{\bb}^{2}
f_{\V0}(\bb_{0})$ that its eigenvalues are non-degenerate.
The following result holds.

\*

\0{\bf Lemma 18.} {\it The Lebesgue measure of the set
$(0,\e_{0})\setminus\EE^{*}$ is bounded proportionally to
some value $G(\e_{0})$, with $G(\e)=o(\e)$.}

\*

\0{\it Proof.} As in the proof of Lemma 12
we start with the first conditions in \equ(4.2).
By setting $\e=\e(t)$, with $t\in[-1,1]$, and defining
$F(t)=\oo\cdot\nn \pm \sqrt{\ul{\l}^{[n]}_{j}(\e(t),\oo)}=
t C_{0} \g^{*}_{n(\nn)}$, one has $|{\rm d}F/{\rm d}t| =
|\dpr_{\e}\sqrt{\ul{\l}^{[n]}_{j}(\e,\oo)}|\,|{\rm d}\e/{\rm d}t| =
C_{0} \g^{*}_{n(\nn)}$, so that, by using that
$\ul{\l}^{[n]}_{j}(\e,\oo) = a_{j}\e + O(\e^{2})$, one finds
$|{\rm d}\e/{\rm d}t| \le B C_{0} \sqrt{\e}\g^{*}_{n(\nn)}$
for some constant $B$. Again, for fixed $\nn$ and $i$,
by the closeness property of the self-energies,
we can impose only the conditions corresponding
to the scales up to $n_{1}(\nn)=O(\log\log 1/\g^{*}_{n(\nn)})$,
at the price of enlarging the sets of excluded values (by a factor 2).
Hence the measure of the set of excluded values $\e\in\EE^{[n-1]}$,
$n \ge n_{0}$, found by imposing the first conditions \equ(4.2)
can be bounded by
$$ \eqalign{
& {\rm const.} \sum_{i=r+1}^{d} C_{0}\sqrt{\e_{0}}\sum_{n=n_{0}}^{\io}
\sum_{n'=n_{0}}^{n_{1}(\nn)}
\sum_{2^{n-1}<|\nn| \le 2^{n}} \g^{*}_{n} \cr
& \qquad \qquad \le
{\rm const.} \, s \, C_{0} \sqrt{\e_{0}}
\sum_{n=n_{0}}^{\io} \g^{*}_{n}\, 2^{nr}
\log \log {1 \over \g^{*}_{n} } \= G (\e_{0}) , \cr}
\Eq(4.9) $$
where we have used that $n_{0}$ is uniquely determined by $\e_{0}$.
Therefore we have
$$ \eqalign{
{ G(\e_{0}) \over \e_{0} } & \le {\rm const.} \,
{1 \over \g_{n_{0}}} \sum_{n=n_{0}}^{\io} \g^{*}_{n}(\oo) \,
2^{nr} \log \log {1 \over \g^{*}_{n} } \cr
& \le {\rm const.} \, {1 \over \g_{n_{0}}(\oo) }
\sum_{n=n_{0}}^{\io} n \, \g^{*}_{n_{0}} 2^{-n(r+1)} 2^{nr} , \cr}
\Eq(4.10) $$
which tends to $0$ as $n_{0}\to \io$ (that is as $\e_{0}\to 0$).

The estimate of the measure of the set of
excluded values $\e\in \EE^{[n-1]}$, $n\ge n_{0}$,
found by imposing the second conditions \equ(4.2) can be obtained
by reasoning in the same way. In such a case we need a lower
bound on $\dpr_{\e} (\sqrt{\ul{\l}^{[n]}_{i}(\e,\oo)}
\pm \sqrt{\ul{\l}^{[n]}_{j}(\e,\oo)})$, which requires $|a_{i}-a_{j}|>a$
for all $i\neq j$ and for some constant $a$;
cf. Ref.~\cita{GG2}, Appendix A2, for a similar discussion. \qed

\*\*
\appendix(A1,Proof of Lemma 10)

\0If we consider two rotation vectors $\oo,\oo'\in\O(C_{0})$, they
are characterised by the respective sequences
$\{\g_{n}(\oo)\}_{n=0}^{\io}$ and $\{\g_{n}(\oo')\}_{n=0}^{\io}$.
For all $n\ge 0$ one has $\g_{n}(\oo) \ge \g_{n}^{*}$ and
$\g_{n}(\oo') \ge \g_{n}^{*}$.

We introduce some shortened notations, by setting
$$ \eqalign{
\X_{n}(\nn,\oo) & = \prod_{p=0}^{n}
\chi_{p}(\D^{[p]}(\oo\cdot\nn;\e,\oo)) , \cr
\Psi_{n}(\nn,\oo) & = \left(
\prod_{p=0}^{n-1} \chi_{p}(\D^{[p]}(\oo\cdot\nn;\e,\oo)) \right)
\psi_{n}(\D^{[n]}(\oo\cdot\nn;\e,\oo)) , \cr}
\Eqa(A1.1) $$
and
$$ \eqalignno{
\X_{n,s}(\nn,\oo',\oo) & =
\Big( \prod_{p=0}^{s-1} \chi_{p}(\D^{[p]}(\oo'\cdot\nn;\e,\oo'))
\Big) \left( \chi_{s}(\D^{[s]}(\oo'\cdot\nn;\e,\oo')) -
\chi_{s}(\D^{[s]}(\oo\cdot\nn;\e,\oo)) \right) \cr
& \qquad
\Big( \prod_{p=s+1}^{n} \chi_{p}(\D^{[p]}(\oo\cdot\nn;\e,\oo))
\Big) , \qquad 0 \le s \le n ,
& \eqa(A1.2) \cr
\Psi_{n,s}(\nn,\oo',\oo) & =
\Big( \prod_{p=0}^{s-1} \chi_{p}(\D^{[p]}(\oo'\cdot\nn;\e,\oo'))
\Big) \left( \chi_{s}(\D^{[s]}(\oo'\cdot\nn;\e,\oo')) -
\chi_{s}(\D^{[s]}(\oo\cdot\nn;\e,\oo)) \right) \cr
& \qquad
\Big( \prod_{p=s+1}^{n-1} \chi_{p}(\D^{[p]}(\oo\cdot\nn;\e,\oo))
\Big) \psi_{n}(\D^{[n]}(\oo\cdot\nn;\e,\oo)) ,
\qquad 0 \le s \le n-1 ,\cr
\Psi_{n,n}(\nn,\oo',\oo) & =
\Big( \prod_{p=0}^{n-1} \chi_{p}(\D^{[p]}(\oo'\cdot\nn;\e,\oo'))
\Big) \left(
\psi_{n}(\D^{[n]}(\oo'\cdot\nn;\e,\oo')) -
\psi_{n}(\D^{[n]}(\oo\cdot\nn;\e,\oo)) \right) , \cr}
$$
where all products have to be meant as 1 when containing no factor.

Finally we define the Hermitian matrices
$$ D_{n}(\nn,\oo) = (\oo\cdot\nn)^{2}-\MM^{[\le n]}(\oo\cdot\nn;\e,\oo) ,
\Eqa(A1.3) $$
with $\MM^{[\le n]}(\oo\cdot\nn;\e,\oo)$ given by \equ(3.3).
In the obvious way one defines also $D_{n}(\nn,\oo')$.

Note that one has
$$ \sum_{p=0}^{n-1} (\g_{p}^{*})^{-m} \le
n \, (\g_{n-1}^{*})^{-m} < (\g_{n}^{*})^{-(m+1)} ,
\Eqa(A1.4) $$
for all $m\in\NNN$.

\*

\0{\bf Lemma A1.} {\it One has
$$ \chi_{n}(x') - \chi_{n}(x) = b(x) \left( x' - x \right) +
o (|x'-x|) , \qquad
\left| b(x) \right| \le \Phi \, (\g^{*}_{n})^{-2} |x'-x| ,
\Eqa(A1.5) $$
for a suitable positive constant $\Phi$.}

\*

\0{\it Proof.} We can write
$$ \left| \chi_{n}(x') - \chi_{n}(x) \right|
\le  \b^{-2} (\g^{*}_{n})^{-2}
\left| x' - x \right| \int_{0}^{1} {\rm d} t \, \dpr_{x} \chi(x(t)) ,
\Eqa(A1.6) $$
where $x(t) = \b^{-2} (\g^{*}_{n})^{-2}(x+ t (x'-x))$
and $|\dpr_{x}\chi(x(t))| \le {\rm const.}$ \qed

\*

By noting that $\psi_{n}=1-\chi_{n}$, we see that Lemma A1 yields the
same bounds as \equ(A1.6) also if we replace $\chi_{n}$ with $\psi_{n}$.

\*

\0{\bf Lemma A2.} {\it For $\oo,\oo'\in\O_{*}(C_{0})$
assume that the bounds \equ(3.28) hold for all $n'\le n$.
Then one has 
$$ \eqalignno{
& \| g^{[n]}(\oo'\cdot\nn;\e,\oo') - g^{[n]}(\oo\cdot\nn;\e,\oo)
- \dpr_{\oo} g^{[n]}(\oo\cdot\nn;\e,\oo)\,(\oo'-\oo) \| =
(\g^{*}_{n})^{-\d} |\nn| \, o(|\oo'-\oo|) , \cr
& \| \dpr_{\oo} g^{[n]}(\oo\cdot\nn;\e,\oo) \|
\le D \, (\g^{*}_{n})^{-\d} |\nn| ,
& \eqa(A1.7) \cr} $$
for suitable positive constants $D$ and $\d$.}

\*

\0{\it Proof.} By using the definition \equ(3.5) we have
$$ \eqalignno{
& g^{[n]}(\oo'\cdot\nn;\e,\oo') - g^{[n]}(\oo\cdot\nn;\e,\oo)
= \Psi_{n}(\nn,\oo') D_{n}^{-1} (\nn,\oo') -
\Psi_{n}(\nn,\oo) D_{n}^{-1} (\nn,\oo) \cr
& \qquad = - \Psi_{n}(\nn,\oo) D_{n}^{-1} (\nn,\oo')
\left( D_{n}(\nn,\oo') - D_{n}(\nn,\oo) \right)
D_{n}^{-1} (\nn,\oo) \cr
& \qquad \qquad + \sum_{p=0}^{n} \Psi_{n,p}(\nn,\oo',\oo)
D^{-1}_{n}(\nn,\oo') ,
& \eqa(A1.8) \cr}
$$
where we can write
$$ \eqalignno{
& D_{n}(\nn,\oo') - D_{n}(\nn,\oo)
& \eqa(A1.9) \cr
& \qquad = (\oo'\cdot\nn+\oo\cdot\nn)
((\oo'-\oo)\cdot\nn) - \MM^{[\le n]} (\oo'\cdot\nn;\e,\oo') +
\MM^{[\le n ]}(\oo\cdot\nn;\e,\oo) \cr}
$$
so that we obtain
$$ \eqalignno{
& D_{n}(\nn,\oo') - D_{n}(\nn,\oo)
& \eqa(A1.10) \cr 
& \qquad = \left( 2 \nn (\oo\cdot\nn) -
\dpr_{\oo}\MM^{[\le n]}(\oo\cdot\nn,\e,\oo) \right) \cdot
\left( \oo'-\oo \right) + ((\oo'-\oo)\cdot\nn)^{2} +
\e^{2} |\nn| \, o(|\oo'-\oo|) , \cr}
$$
by the assumed estimate \equ(3.28).

If $\Psi_{n}(\nn,\oo') \neq 0$ we can bound
the last sum in \equ(A1.8) by
$$ \eqalignno{
\sum_{p=0}^{n} \Psi_{n,p}(\nn,\oo',\oo)
\left\| D^{-1}_{n}(\nn,\oo') \right\|
& \le {\rm const.} \, (\g^{*}_{n})^{-2} \sum_{p=0}^{n} 
(\g^{*}_{p})^{-2} \, |\nn| \, |\oo'-\oo| \cr
& \le {\rm const.} \, (\g^{*}_{n})^{-5} |\nn| \, |\oo'-\oo| ,
& \eqa(A1.11) \cr} $$
where we have used Lemma \secc(A1) to bound $\Psi_{n,p}(\nn,\oo',\oo)$,
and \equ(A1.4) to perform the sum over $p=0,\ldots,n$.
Note that in order to profitably use the bound \equ(A1.5) we have
to use that the bounds \equ(3.28) and the consequent \equ(3.29) imply
analogous bounds also for the propagator divisors $\D^{[n]}(x;\e,\oo)$,
without the factor $\e^{2}$: this follows from the fact that
\equ(3.4) defines functions which are smooth in $\e$ and $\oo$.

Still if $\Psi_{n}(\nn,\oo') \neq 0$ we can bound
in \equ(A1.8) also the matrices $D_{n}^{-1} (\nn,\oo')$ and
$D_{n}^{-1} (\nn,\oo)$ both proportionally to
$(\g^{*}_{n})^{-2}$, so that \equ(A1.7) follows.

If $\Psi_{n}(\nn,\oo') = 0$ call $\a_{n}$ and $\a_{n}'$ the
eigenvalues with minimum absolute value of
$D_{n}(\nn,\oo)$ and $D_{n}(\nn,\oo')$, respectively.
If $|\a_{n}'|\ge |\a_{n}|$ we can proceed as in the previous case,
and we obtain the same bound.

Finally if $\Psi_{n}(\nn,\oo') = 0$ and $|\a_{n}'|< |\a_{n}|$,
we can write
$$ g^{[n]}(\oo'\cdot\nn;\e,\oo') - g^{[n]}(\oo\cdot\nn;\e,\oo) =
- g^{[n]}(\oo\cdot\nn;\e,\oo) .
\Eqa(A1.12) $$
Moreover, by the assumed bound \equ(3.28) we have that the
difference $D_{n}(\nn,\oo') - D_{n}(\nn,\oo)$
is given as in \equ(A1.10), so that
$\| D_{n}(\nn,\oo') - D_{n}(\nn,\oo) \| \le 4|\nn|\,|\oo'-\oo|$.
Hence the difference between the eigenvalues of 
$D_{n}(\nn,\oo)$ and $D_{n}(\nn,\oo')$ is bounded by
$C|\nn|\,|\oo'-\oo|$, for some constant $C$; this is again
a consequence of Lidski\u\i's theorem.
Therefore for $|\a_{n}| \ge 2 C |\nn| \, |\oo'-\oo|$ we can bound
$|\a_{n}'| \ge |\a_{n}| - C|\nn|\,|\oo-\oo'| \ge |\a_{n}|/2$, and
\equ(A1.7) follows once more by reasoning as before, whereas for
$|\a_{n}| < 2 C |\nn| \, |\oo'-\oo|$ we can bound in \equ(A1.12)
$$ \eqalign{
\left\| g^{[n]}(\oo\cdot\nn;\e,\oo) \right\| & \le
\left\| D_{n}^{-1}(\nn,\oo) \right\| \le
{\left\| D_{n}^{-1}(\nn,\oo) \right\| \over |\a_{n}| } |\a_{n}| \cr
& \le {1 \over |\a_{n}|^{2}} 2C|\nn|\,|\oo'-\oo| \le
{\rm const.} (\g_{n}^{*})^{-2} |\nn|\,|\oo'-\oo| , \cr}
\Eqa(A1.13) $$
and \equ(A1.7) follows also in such a case. \qed

\*

Now we can prove Lemma 10. The proof is by induction on $n$.
One can reason as in Ref.~\cita{Ge3}. More precisely one writes
$$ \eqalignno{
& \MM^{[\le n]}(\oo'\cdot\nn;\e,\oo') -
\MM^{[\le n]}(\oo\cdot\nn;\e,\oo)
& \eqa(A1.14) \cr
& = \sum_{p=0}^{n} \X_{p}(\nn,\oo)
\left( M^{[p]}(\oo'\cdot\nn;\e,\oo') -
M^{[p]}(\oo\cdot\nn;\e,\oo) \right) \cr
& \qquad \qquad +
\sum_{p=0}^{\n} \sum_{s=0}^{p}
\X_{p,s}(\nn,\oo',\oo) \, M^{[p]}(\oo'\cdot\nn;\e,\oo') , \cr}
$$
where the matrices $M^{[p]}(x;\e,\oo)$ are defined
as in Ref.~\cita{GG2}, formula (5.9). In the first term of the sum
the difference $M^{[p]}(\oo'\cdot\nn;\e,\oo') -
M^{[p]}(\oo\cdot\nn;\e,\oo)$ can be written
as a sum over self-energy clusters $T$ of differences
of self-energy values $\VV_{T}(\oo'\cdot\nn)-\VV_{T}(\oo\cdot\nn)$,
computed with rotation vectors $\oo'$ and $\oo$, respectively.
The latter can be written as a sum of several contributions,
each of which is given by the product of a factor $\AAA(\oo')$
depending on $\oo'$ but no on $\oo$ times a
factor\footnote{${}^2$}{\nota No relation with the Bryuno
function $\BBB(\o)$ introduced in Section \secc(2).} $\BBB(\oo)$
depending on $\oo$ but not on $\oo'$ times a difference of propagators
$\D_{\ell}(\oo',\oo)= g^{[n_{\ell}]}(\oo'\cdot\nn;\e,\oo')-
g^{[n_{\ell}]}(\oo\cdot\nn;\e,\oo)$, with $n_{\ell}<n$.

The difference $\D_{\ell}(\oo',\oo)$ can be bounded
according to Lemma A4, proportionally to $|\nn| \, |\oo'-\oo|$,
by using that $n_{\ell}<n$ and the inductive hypothesis.

Moreover the decomposition $\AAA(\oo')\D_{\ell}(\oo',\oo)\BBB(\oo)$
can be made in such a way that the factor $\AAA(\oo')$ corresponds
to a connected subset $T_{0}$ of $T$, while the factor $\BBB(\oo)$
is the product of factorising factors $\BBB_{i}(\oo)$ corresponding to
subsets $T_{i}$ of $T$ containing lines preceding the lines of $T_{0}$
(again we refer to Ref.~\cita{GG1} and \cita{Ge3} for details).
The only factor which requires some care is that corresponding
to the subset, say $T_{1}$, connected to the entering line of $T$.
But this can be easily discussed as in deriving Lemma 4.
Indeed one can prove by induction (on the number of nodes) that,
given a subset $T_{1}$ with the considered structure, one has
$N_{p}(T_{1}) \le K \, 2^{-p}M(T_{1})$ for all $p<n$ (note the
absence of the summand $-1$ with the respect the analogous
inductive assumption one makes for $N_{h}(\th)$ and $N_{h}(T)$).

Then, by taking into account also the factor $(\g^{*}_{n_{\ell}})^{-\d}$
possibly arising from $\D_{\ell}(\oo',\oo)$, each factor $\BBB_{i}(\oo)$
can be bounded by $B_{1}^{k_{i}}{\rm e}^{-k|\nn_{i}|/2}$,
where $k_{i}$ is the number of vertices in $T_{i}$, $\nn_{i}$
is the momentum of the line $\ell_{i}$ connecting $T_{i}$ to $T_{0}$,
and $B_{1}$ is some positive constant (the proof proceeds as for Lemma 9).
The factor $\AAA(\oo')$ can be bounded in the same way, and,
by taking into account also the factor $(\g^{*}_{n_{\ell}})^{-\d}$
possibly arising from $\D_{\ell}(\oo',\oo)$
and the factors ${\rm e}^{-k|\nn_{i}|/2}$, it can be bounded
by $B_{2}^{k_{0}}$, where $k_{0}$ is the number of vertices in $T_{0}$,
and $B_{2}$ is some other positive constant (again the proof
proceeds as for Lemma 9, but with $\k$ replaced with $\k/2$).
By writing $\AAA(\oo')=\AAA(\oo)+(\AAA(\oo')-\AAA(\oo))$
one can iterate the construction above for the difference
$\AAA(\oo')-\AAA(\oo)$. The only difference with respect to the previous
case is that now the factor $\k/2$ is replaced with $\k/4$.

All the other terms of the double sum in \equ(A1.14) can be discussed
in a similar way, by relying once more on Lemma \secc(A2).
We omit the details, which can be worked out as in Ref.~\cita{Ge3}.

Therefore the property \equ(3.28) of Lemma 10 is proved.

\*\*
\appendix(A2,Extensions and generalisations)

\0The analysis performed in this article applies
to more general Hamiltonians of the form
$$ \HHH = \HHH_{0}(\III) + \e f(\III,\ff) ,
\Eqa(A2.1) $$
where $(\III,\ff) \in \AAAA \times \TTT^{d}$,
with $\AAAA \subset \RRR^{d}$, are conjugate action-angle variables,
$\e$ is a real parameter, and the functions $\HHH_{0}$ and $f$
are assumed to be analytic in their arguments.
We assume also convexity on $\HHH_{0}$ and a non-degeneracy
condition on $f$ which will be specified later.
Here we confine ourselves to sketch the basic arguments:
full details will be published elsewhere \cita{GGG2}.
Note that in principle weaker conditions on the unperturbed
Hamiltonian could be considered; we refer to Ref.~\cita{Se}
for a survey of results under the usual Diophantine conditions.

We say that a vector $\oo_{*}\in\RRR^{d}$ is $s$-resonant
if it satisfies $s$ resonance conditions
$\nn_{i}\cdot\oo_{*}=0$, for $s$ linearly independent integer
vectors $\nn_{1},\ldots,\nn_{s}$.

With respect to \equ(A2.1), the Hamiltonian \equ(1.1) has a very
special form. Even by considering Hamiltonians of the form
$$ \HHH = {1\over2}\III\cdot\III + \e f(\ff) ,
\Eqa(A2.2) $$
it is not always possible to reduce to \equ(1.1).
In fact if $\II \! : \RRR^{d} \to \RRR^{d}$
is the linear operator which transforms $\oo_{*}$ into $(\oo,\V0) \in
\RRR^{r} \times \RRR^{s}$, where $\oo$ has rationally independent
components, then the action variables $(\AA,\BB)$
are mixed together and also terms of the form $A_{i} B_{j}$ appear.

Though, it is easy to extend the analysis to such a case.
And with a little further work, we can consider also Hamiltonians with
any unperturbed Hamiltonian $\HHH_{0}$ satisfying a convexity property
(so that the eigenvalues of the matrix $\det \dpr_{\III}^{2}\HHH_{0}$
are all strictly positive). The frequency map $\III\to \oo =
\dpr \HHH_{0}/\dpr\III$ is a local diffeomorphsism, so that if we
fix $\III_{0}$ in such a way that the corresponding rotation
vector $\oo(\III)$ is $s$-resonant, we can find an immersed
$r$-dimensional manifold $\MMMM_{r}$, with $r=d-s$,
containing $\III_{0}$, on which the $s$ resonance conditions
are satisfied. We shall call $\MMMM_{r}$ a {\it resonant manifold} 

Under the action of the symplectic transformation given by the
lift $\SS$ of $\II$, we can pass to new coordinates, which we continue
to denote with the same symbols, such that in the new coordinates
the rotation vector has became $(\oo,\V0)$. For simplicity we still
call $\MMMM_{r}$ the resonant manifold in the new coordinates.

As we are interested in local properties (in the action variables)
we can assume that a system of coordinates adapted to $\MMMM_{r}$
has been fixed, so that we can write $\III=(\AA,\BB)$ in such a way that
$\BB=\V0$ identifies $\MMMM_{r}$. For $\e=0$ a motion
on $\MMMM_{r}$ is determined by fixing $\BB=\V0$ and
$\AA=\AA_{0}$ in such a way that the conjugated angles $(\aa,\bb)$
move according to the law $(\aa,\bb) \to (\aa+\oo t,\bb)$,
with $\oo$ uniquely determined by $\AA_{0}$.
This means that the unperturbed lower-dimensional
tori can be characterized by the rotation vectors $\oo\in\RRR^{r}$
depending on the action variables $\AA$. Hence for $\e\neq 0$
the Hamiltonian describing the system can be written as
$$ \HHH = \HHH_{0}(\AA,\BB) + \e f(\AA,\BB,\aa,\bb) ,
\Eqa(A2.3) $$
so that the subset of $\MMMM_{r}\times\TTT^{r}$ whose
unperturbed invariant tori can be continued under
perturbation can be characterised by the set of allowed values
of $\oo$, provided the map $\AA \mapsto \oo(\AA)$
is a local diffeomorphism, which is true under our hypotheses.

When the perturbation depends also on the action variables,
as in \equ(A2.3), of course one needs both equations
for action and angle variables:
$$ \cases{
\dot \AA = - \e \dpr_{\aa}f(\AA,\BB,\aa,\bb) , & \cr
\dot \BB = - \e \dpr_{\bb}f(\AA,\BB,\aa,\bb) , & \cr
\dot \aa = \dpr_{\AA}\HHH_{0}(\AA,\BB) +
\e \dpr_{\AA}f(\AA,\BB,\aa,\bb) , & \cr
\dot \bb = \dpr_{\BB}\HHH_{0}(\AA,\BB) +
\e \dpr_{\BB}f(\AA,\BB,\aa,\bb) . & \cr}
\Eqa(A2.4) $$
The main difference with respect to the
analysis in Sections \secc(3) and \secc(4) is that
the propagators are of the form \equ(3.4), with
$(x^{2}-\MM^{[\le n]}(x;\e,\o))^{-1}$ replaced with
$(ix-\MM^{[\le n]}(x;\e,\o))^{-1}$, where we can write
$\MM^{[\le n]}(x;\e,\o)=\CC+\NN^{[\le n]}(x;\e,\o)$,
if the matrix $\CC$ is given by
$$ \CC = \left( \matrix{
0 & 0 & 0 & 0 \cr
0 & 0 & 0 & 0 \cr
\dpr_{\AA}^{2}\HHH_{0} & \dpr_{\AA}\dpr_{\BB}\HHH_{0} & 0 & 0 \cr
\dpr_{\AA}\dpr_{\BB}\HHH_{0} & \dpr_{\BB}^{2}\HHH_{0} & 0 & 0 \cr}
\right) ,
\Eqa(A2.5) $$
and the matrix $\NN^{[\le n]}\=\NN^{[\le n]}(x;\e,\o)$ is such that,
by extracting the dominant order, and setting $x=0$, one has
$$ \eqalign{
\NN^{[\le n]}(0;\e,\oo) & =
\NN^{[0]}(0;\e,\oo) + O(\e^{2}) , \cr
\NN^{[0]}(0;\e,\oo) & = \left( \matrix{
0 & 0 & 0 & 0 \cr
- \e \dpr_{\bb}\dpr_{\AA} f &
- \e \dpr_{\bb}\dpr_{\BB} f & 0 &
- \e \dpr_{\bb}^{2} f \cr
\e \dpr_{\AA}^{2} f &
\e \dpr_{\AA}\dpr_{\BB} f & 0 &
\e \dpr_{\AA}\dpr_{\bb} f \cr
\e \dpr_{\BB}\dpr_{\AA} f &
\e \dpr_{\BB}^{2} f &  0 &
\e \dpr_{\BB}\dpr_{\bb} f \cr} \right) , \cr}
\Eqa(A2.6) $$
whereas the other terms depending on $x$ which are not negligible
with respect with the dominant ones appear as
$$ \NN^{[\le n]}_{1}(x;\e,\oo) \= \left( \matrix{
O(\e^{2} x) & O(\e^{2} x) & O(\e^{2} x^{2}) & O(\e^{2} x) \cr
0 & 0 & O(\e^{2} x) & 0 \cr
0 & 0 & 0(\e^{2} x) & 0 \cr
0 & 0 & 0(\e^{2} x) & 0 \cr} \right) .
\Eqa(A2.7) $$
Some deep relations turn out to exist between the
matrices $\MM^{[\le n]}(x;\e,\oo)\,E$ and their transposed,
if $E$ denotes the standard symplectic matrix.
Then, by using these relations, one can bound the propagators
in terms of the eigenvalues
of a suitable symplectic matrix $S$: for the latter,
besides $d$ harmless eigenvalues of order 1 (in $\e$ and $x$)
there are $r$ eigenvalues proportional to $x^{2}$,
while the other $s$ eigenvalues are of the form
$x^{2} - \l^{[0]}_{j}(x,\e,\oo) + O(\e x) + O(\e^{2})$ , with
$\l^{[0]}_{j}(x,\e,\oo)= \e a_{j-r}(\oo)$, $j=r+1,\ldots,d$,
if $a_{1}(\oo),\ldots,a_{s}(\oo)$ are the dominant terms
of the normal frequencies (which depend also on $\AA$, hence
on $\oo$, in this case). The aforementioned non-degeneracy condition
on $f$ is that the functions $a_{j}(\oo)$ are all strictly positive.
The dependence on $\oo$ does not introduce any further difficulties, and
in fact Whitney differentiability in $\oo$ (as it appears in the
subsequent iterative steps) would be enough to carry on the analysis.

Hence, the situation is very similar to that which has been
considered in the previous Sections \secc(3) and \secc(4),
and one can reason essentially as there.
Of course notations become more cumbersome, because
one has to keep trace also of the action variables
(which cannot be any more expressed trivially in terms
of the angle variables), and more sophisticated diagrammatic
rules have to be envisaged; again see Ref.~\cita{GGG2}
for details. But the basic estimates and arguments
are the same, and the same conclusions hold.

\*\*
\0{\bf Acknowledgments.} I am indebted to A. Giuliani for
many discussions and comments.

\*\*

\rife{BaG}{1}{
M. Bartuccelli, G. Gentile,
{\it Lindstedt series for perturbations of isochronous systems.
A review of the general theory},
Rev. Math. Phys.
{\bf 14} (2002), no. 2, 121--171. }
\rife{BeG}{2}{
A. Berretti, G. Gentile,
{\it Bryuno function and the standard map},
Comm. Math. Phys.
{\bf 220} (2001), no. 3, 623--656. }
\rife{Bo1}{3}{
J. Bourgain,
{\it Construction of quasi-periodic solutions for Hamiltonian
perturbations of linear equations and applications to nonlinear PDE},
Internat. Math. Res. Notices
{\bf 11} (1994), 475--497. }
\rife{Bo2}{4}{
J. Bourgain,
{\it On Melnikov's persistency problem},
Math. Rev. Lett.
{\bf 4} (1997), no. 4, 445--458. } 
\rife{BKS}{5}{
J. Bricmont, A. Kupiainen, A. Schenkel,
{\it Renormalization group and the Melnikov problem for PDE's},
Comm. Math. Phys.
{\bf 221} (2001), no. 1, 101--140. }
\rife{Br1}{6}{
A.D. Bryuno,
{\it Analytic form of differential equations. I},
Trudy Moskov. Mat. Ob\v s\v c.
{\bf 25} (1971), 119--262;
english translation in
Trans. Moscow Math. Soc.
{\bf 25} (1973), 131--288. }
\rife{Br2}{7}{
A.D. Bryuno,
{\it Analytic form of differential equations. II},
Trudy Moskov. Mat. Ob\v s\v c.
{\bf 26} (1972), 199--239;
english translation in
Trans. Moscow Math. Soc. 
{\bf 26} (1974), 199--239. }
\rife{Ch1}{8}{
Ch.-Q. Cheng,
{\it Birkhoff-Kolmogorov-Arnold-Moser tori in convex Hamiltonian systems},
Comm. Math. Phys.
{\bf 177} (1996), no. 3, 529--559. }
\rife{Ch2}{9}{
Ch.-Q. Cheng,
{\it Lower-dimensional invariant tori in the regions of 
instability for nearly integrable Hamiltonian systems},
Comm. Math. Phys.
{\bf 203} (1999), no. 2, 385--419. }
\rife{ChW}{10}{
Ch.-Q. Cheng, Sh. Wang,
{\it The surviving of lower-dimensional tori from a
resonant torus of Hamiltonian systems},
J. Differential Equations
{\bf 155} (1999), no. 2, 311--326. }
\rife{CG}{11}{
L. Chierchia, G. Gallavotti
{\it Smooth prime integrals for quasi-integrable Hamiltonian systems},
Nuovo Cimento B (11)
{\bf 67} (1982), no. 2, 277--295. } 
\rife{D}{12}{
A.M. Davie,
{\it The critical function for the semistandard map},
Nonlinearity 
{\bf 7} (1994), no. 1, 219--229. }
\rife{EV}{13}{
J. Ecalle, B. Vallet,
{\it Bruno Correction and linearization of resonant
vector fields and diffeomorphisms},
Math. Z.
{\bf 229} (1998), no. 2, 249--318. }
\rife{E1}{14}{
L.H. Eliasson,
{\it Perturbations of stable invariant tori for Hamiltonian systems},
Ann. Scuola Norm. Sup. Pisa Cl. Sci. (4)
{\bf 15} (1988), no. 1, 115--147. }
\rife{E2}{15}{
L.H. Eliasson,
{\it Hamiltonian systems with linear normal form near an invariant torus},
Nonlinear dynamics (Bologna, 1988),  11--29,
World Sci. Publishing, Teaneck, NJ, 1989. }
\rife{Ga1}{16}{
G. Gallavotti,
{\it  A criterion of integrability for perturbed nonresonant
harmonic oscillators. ``Wick ordering" of the perturbations
in classical mechanics and invariance of the frequency spectrum},
Comm. Math. Phys.
{\bf 87} (1982/83), no. 3, 365--383. }
\rife{GG1}{17}{
G. Gallavotti, G. Gentile,
{\it Hyperbolic low-dimensional tori and summations of divergent series},
Comm. Math. Phys.
{\bf 227} (2002), no. 3, 421--460. }
\rife{GGG1}{18}{
G. Gallavotti, G. Gentile, A. Giuliani,
{\it Fractional Lindstedt series},
Preprint, 2005. }
\rife{GGG2}{19}{
G. Gallavotti, G. Gentile, A. Giuliani,
work in preparation. }
\rife{Ge1}{20}{
G. Gentile,
{\it Diagrammatic techniques in perturbations theory, and applications},
Symmetry and perturbation theory (Rome, 1998),  59--78,
World Sci. Publishing, River Edge, NJ, 1999. }
\rife{Ge2}{21}{
G. Gentile,
{\it Quasi-periodic solutions for two-level systems},
Comm. Math. Phys.
{\bf 242} (2003), no. 1-2, 221-250. }
\rife{Ge3}{22}{
G. Gentile,
{\it Pure point spectrum for two-level systems in a strong
quasi-periodic field},
J. Statist. Phys.
{\bf 115} (2004), no. 5-6, 1605-2620. }
\rife{GG2}{23}{
G. Gentile, G. Gallavotti,
{\it Degenerate elliptic tori},
Comm. Math. Phys.
{\bf 257} (2005), no. 2, 319--362. }
\rife{GM}{24}{
G. Gentile, V. Mastropietro,
{\it Methods for the analysis of the Lindstedt series for KAM tori
and renormalizability in classical mechanics. A review with
some applications},
Rev. Math. Phys.
{\bf 8} (1996), no. 3, 393--444. }
\rife{GL}{25}{
A. Giorgilli, U. Locatelli,
{\it On classical series expansions for quasi-periodic motions},
Math. Phys. Electron. J.
{\bf 3} (1997), Paper 5, 25 pp. }
\rife{JLZ}{26}{
\`A. Jorba, R. de la Llave, M. Zou,
{\it Lindstedt series for lower-dimensional tori}, in
{\it Hamiltonian systems with three or more degrees of
freedom} (S'Agar\'o, 1995), 151--167,
NATO Adv. Sci. Inst. Ser. C Math. Phys. Sci., 533, Ed. C. Sim\'o,
Kluwer Acad. Publ., Dordrecht, 1999. }
\rife{JV1}{27}{
\`A. Jorba, J. Villanueva,
{\it On the persistence of lower-dimensional invariant tori
under quasi-periodic perturbations},
J. Nonlinear Sci.
{\bf 7} (1997), no. 5, 427--473. }
\rife{JV2}{28}{
\`A. Jorba, J. Villanueva,
{\it On the normal behaviour of partially elliptic
lower-dimensional tori of Hamiltonian systems},
Nonlinearity
{\bf 10} (1997), no. 4, 783--822. }
\rife{Ka}{29}{
T. Kato,
{\it Perturbation theory for linear operators},
Classics in Mathematics, Springer, Berlin, 1995. }
\rife{Ko}{30}{
A.N. Kolmogorov,
{\it On conservation of conditionally periodic motions
for a small change in Hamilton's function},
Dokl. Akad. Nauk SSSR
{\bf 98} (1954). 527--530. }
\rife{Ku1}{31}{
S.B. Kuksin,
{\it Hamiltonian perturbations of infinite-dimensional linear systems
with imaginary spectrum},
Akade\-miya Nauk SSSR. Funktsional$^{\prime}$ny\u\i\ Analiz i ego Prilozheniya
{\bf 21} (1987),  no. 3, 22--37. }
\rife{Ku2}{32}{
S.B. Kuksin,
{\it Nearly integrable infinite-dimensional Hamiltonian systems},
Lecture Notes in Mathematics, 1556, Springer-Verlag, Berlin, 1993. }
\rife{L}{33}{
V.B. Lidski\u\i,
{\it On the characteristic numbers of the sum and
product of symmetric matrices},
Dokl. Akad. Nauk SSSR (N.S.)
{\bf 75} (1950), 769--772. }
\rife{LD}{34}{
J. Lopes Dias,
{\it Renormalization and reducibility of Bryuno skew-systems},
Preprint, 2004. }
\rife{M1}{35}{
V.K. Mel$'$nikov,
{\it On certain cases of conservation of conditionally periodic motions
under a small change of the Hamiltonian function},
Dokl. Akad. Nauk SSSR
{\bf 165} (1965), 1245--1248;
english translation in
Soviet Math. Dokl.
{\bf 6} (1965), 1592--1596. }
\rife{M2}{36}{
V.K. Mel$'$nikov,
{\it A certain family of conditionally periodic solutions of a
Hamiltonian systems},
Dokl. Akad. Nauk SSSR
{\bf 181} (1968), 546--549;
english translation in
Soviet Math. Dokl.
{\bf 9} (1968), 882-886. }
\rife{P1}{37}{
J. P\"{o}schel,
{\it Integrability of Hamiltonian systems on Cantor sets},
Comm. Pure Appl. Math.
{\bf 35} (1982), no. 5, 653--696. }
\rife{P2}{38}{
J. P\"{o}schel,
{\it On elliptic lower-dimensional tori in Hamiltonian systems},
Math. Z.
{\bf 202} (1989), no. 4, 559--608. }
\rife{Re}{39}{
F. Rellich,
{\it Perturbation theory of eigenvalue problems},
Gordon and Breach Science Publishers, New York-London-Paris, 1969. }
\rife{R1}{40}{
H. R\"ussmann,
{\it On the one-dimensional Schr\"odinger equation
with a quasiperiodic potential},
Nonlinear dynamics (Internat. Conf., New York, 1979), pp. 90--107,
Ann. New York Acad. Sci., 357,
New York Acad. Sci., New York, 1980. }
\rife{R2}{41}{
H. R\"ussmann,
{\it Invariant tori in non-degenerate nearly
integrable Hamiltonian systems}, 
Regul. Chaotic Dynam.
{\bf 6} (2001), 119--204. }
\rife{S}{42}{
W.M. Schmidt,
{\it Diophantine approximation},
Lecture Notes in Mathematics 785, Springer, Ber\-lin, 1980. }
\rife{Se}{43}{
M.B. Sevryuk,
{\it The classical KAM theory at the dawn of the twenty-first century},
Mosc. Math. J.
{\bf 3} (2003), no. 3, 1113--1144. }
\rife{Th}{44}{
W. Thirring,
{\it A course in mathematical physics.
Vol. I. Classical dynamical systems},
Springer, New York-Vienna, 1978. }
\rife{T}{45}{
D.V. Treshch\"ev,
{\it A mechanism for the destruction of resonance tori
in Hamiltonian systems},
Mat. Sb.
{\bf 180} (1989), no. 10, 1325--1346;
english translation in
Math. USSR-Sb.
{\bf 68} (1991), no. 1, 181--203. }
\rife{W}{46}{
H. Whitney,
{\it Analytic extensions of differential
functions defined in closed sets},
Trans. Amer. Math. Soc.
{\bf 36} (1934), no. 1, 63--89. }
\rife{Xu}{47} {
J. Xu,
{\it Persistence of elliptic lower-dimensional invariant
tori for small perturbation of degenerate integrable Hamiltonian systems},
J. Math. Anal. Appl.
{\bf 208} (1997),  no. 2, 372--387. }
\rife{XY}{48} {
J. Xu, J. You,
{\it Persistence of lower dimensional 
tori under the first Melnikov's non-resonance condition},
J. Math. Pures Appl. (9)
{\bf 80} (2001), no. 10, 1045--1067. }
\rife{Y}{49}{
J.-Ch. Yoccoz,
{\it Th\'eor\`eme de Siegel, nombres de Bruno et polin\^omes quadratiques},
Ast\'e\-risque
{\bf 231} (1995), 3--88. }
\rife{You}{50}{
J. You,
{\it Perturbations of lower-dimensional tori for Hamiltonian systems},
J. Differential Equations
{\bf 152} (1999), no. 1, 1--29. }

\biblio

\end